\documentclass[preprint,12pt,authoryear]{elsarticle}

\usepackage[hidelinks]{hyperref}
\usepackage{amsmath,amsthm,color}
\usepackage{amssymb}  
\usepackage{amsfonts}
\usepackage{graphicx}
\usepackage{mathtools}
\usepackage{dsfont}
\usepackage{psfrag}
\usepackage{subfigure,epstopdf}
\usepackage{multirow}
\usepackage{breqn}
\usepackage{float}
\usepackage{cite}
\usepackage{bbm}
\usepackage{algorithm}
\usepackage{algorithmic}
\usepackage{savesym}
\usepackage{float}
\usepackage{tikz}
\usetikzlibrary{positioning,chains,fit,shapes,calc}
\usetikzlibrary{backgrounds}
\usetikzlibrary{arrows.meta}
\usetikzlibrary{shapes,arrows}
\usepackage{tabularx}
\usepackage{optidef}
\usepackage[font=small]{caption}

\theoremstyle{definition}

\newtheorem{remark}{Remark}

\theoremstyle{plain}

\DeclareMathOperator*{\argmin}{argmin}

\journal{}
\begin{document}

\begin{frontmatter}



\title{Fractional cyber-neural systems -- a brief survey\tnoteref{label_t}}
 
 \tnotetext[label_t]{
 This work was supported in part by FCT project POCI-01-0145-FEDER-031411-HARMONY, National Science Foundation GRFP DGE-1842487, Career Award CPS/CNS-1453860, CCF-1837131, MCB-1936775, CNS-1932620, CMMI-1936624, CMMI-1936578, the University of Southern California Annenberg Fellowship, USC WiSE Top-Off Fellowship, the DARPA Young Faculty Award and DARPA Director Award N66001-17-1-4044. The views, opinions, and/or findings contained in this article are those of the authors and should not be interpreted as representing the official views or policies, either expressed or implied by the Defense Advanced Research Projects Agency, the Department of Defense, or the National Science Foundation.\\$^\ast$The first three authors contributed equally.}


\author[label1]{Emily Reed$^\ast$}
\ead{emilyree@usc.edu}

\author[label2]{Sarthak Chatterjee$^\ast$}
\ead{sarthak.chatterjee92@gmail.com}

\author[label3]{Guilherme Ramos$^\ast$}
\ead{guilhermeramos21@gmail.com}

\author[label1]{Paul Bogdan}
\ead{pbogdan@usc.edu}

\author[label4]{S\'ergio Pequito}
\ead{Sergio.Pequito@tudelft.nl}

\affiliation[label1]{organization={Ming Hsieh Department of Electrical and Computer Engineering, University of Southern California},city={Los Angeles},state={CA},country={USA}}

\affiliation[label2]{organization={Department of Electrical, Computer, and Systems Engineering, Rensselaer Polytechnic Institute},city={Troy},state={NY},country={USA}}  
            
\affiliation[label3]{organization={Department of Electrical and Computer Engineering, Faculty of Engineering, University of Porto},country={Portugal}}
            
\affiliation[label4]{organization={Delft Center for Systems and Control, Delft University of Technology},city={Delft},country={The Netherlands}}

\begin{abstract}
Neurotechnology has made great strides in the last 20 years. However, we still have a long way to go to commercialize many of these technologies as we lack a unified framework to study cyber-neural systems (CNS) that bring the hardware, software, and the neural system together. Dynamical systems play a key role in developing these technologies as they capture different aspects of the brain and provide insight into their function. Converging evidence suggests that fractional-order dynamical systems are advantageous in modeling neural systems because of their compact representation and accuracy in capturing the long-range memory exhibited in neural behavior. In this brief survey, we provide an overview of fractional CNS that entails fractional-order systems in the context of CNS. In particular, we introduce basic definitions required for the analysis and synthesis of fractional CNS, encompassing system identification, state estimation, and closed-loop control. Additionally, we provide an illustration of some applications in the context of CNS and draw some possible future research directions. Ultimately, advancements in these three areas will be critical in developing the next generation of CNS, which will, ultimately, improve people’s quality of life.
\end{abstract}



\begin{keyword}
fractional-order systems \sep cyber-neural systems \sep neurotechnology



\end{keyword}

\end{frontmatter}



{
\scriptsize
\tableofcontents
}
\section{Introduction}

We have witnessed an increase in the popularity of neurotechnology, in part propelled by several Silicon Valley companies such as NeuraLink~\citep{neuralink} (founded by Elon Musk), Google, and Facebook, just to mention a few. This tendency is now emerging in Europe as well with a variety of start-up companies across different countries. Yet, we have a long path moving forward to commercialize these devices to a clinical domain~\citep{lewis_2020,whitepaper,standards}. Among the different neurotechnologies, the one experiencing the biggest growth is that of
neurostimulation devices to assess the neural activity (e.g., by tracking the change in electrical potential) and to inject a timely stimulus (e.g., current in electrical neurostimulation devices) that aims to disrupt such activity~\citep{neuralinkhype}. These devices consist of tightly integrated hardware/software components and together with the neural system that they monitor and regulate they form a cyber-neural system~(CNS).

Neural systems like the brain generally exhibit quite diverse activity patterns across subjects under different operational setups~\citep{bargmann2014brain}. Therefore, it is important to develop and translate the tools to enable CNS to become tomorrow's reality.  Yet, to several scientists and engineers, this is \emph{the 21st century equivalent of the space race that led the man to the moon}. Consequently, it is imperative to establish a unified robust framework to understand and regulate brain activity across individuals and regimes (both healthy as well as diseased/disordered)~\citep{markram2012human, van2013wu,ledoux1998emotional}.  Ultimately, this will enable the improvement of people's quality of life.


Fortunately, each year, we get new insights and understanding about life-changing neurological diseases. These advances made in understanding neural systems that provide adequate treatment for these diseases have been mostly achieved with the help of technology that can measure and record neural activity \citep{fairclough2020grand}. Scientists and researchers use these measurements of the brain's activity to create models of the brain.

There are different methods to perform analysis and design of neural dynamical systems. One useful tool in modeling is  dynamical network systems \citep{bassett2017network}.  For example, \citep{presigny2021multiscale} provides an overview of the recent advances in modeling the multi-scale behavior of the brain using dynamical networks. These models have allowed researchers to draw conclusions regarding the brain's topology and function. While many studies have made use of linear dynamical system models \citep{ashourvan2020model, li2017fragility,pequito2017spectral}, these models are unable to capture the nonlinear and non-Markovian behavior exhibited in the brain \citep{zhang2012remarks,west2016fractional,shlesinger1993strange}. On the other hand, several studies make use of more complex nonlinear models; however, these models are not easy to interpret and explain in the context of brain dynamics \citep{west2015fractional,bonilla2007fractional}. 

Fractional-order dynamical systems originated in physics and economics and quickly found their way into engineering applications. Their appeal is mainly due to their representation as a compact  spatiotemporal dynamical system with two  easy-to-interpret sets of parameters, namely the \emph{fractional-order exponents} and the \emph{spatial coupling}. The fractional-order exponents capture the long-range memory in the dynamics of each state variable of the system and the spatial matrix represents the spatial coupling between different state variables, the latter being described by a matrix.

Fractional-order systems provide an efficient way to model many different systems~\citep{valerio2013fractional,west2014colloquium,kilbas2006theory,baleanu2012fractional,podlubny1998fractional,sabatier2007advances,xue2017constructing,lundstrom2008fractional,werner2010fractals,turcott1996fractal,thurner2003scaling,teich1997fractal,chen2010anomalous,petravs2011fractional}. For example, fractional-order systems have been used in domains as disparate as biological networks~\citep{west2016networks}, cyber-physical systems~\citep{xue2016spatio}, nanotechnology~\citep{baleanu2010new}, finance~\citep{scalas2000fractional}, quantum mechanics~\citep{shahin2009fractional}, phasor measurement unit~(PMU) data in the power grid~\citep{shalalfeh2020fractional}, and networked control systems~\citep{cao2009distributed,chen2010fractional,ren2011distributed}, to mention a few.

In this brief survey, we focus our attention on neural behavior, which can be accurately represented by fractional-order systems~\citep{baleanu2011fractional,west2016networks,moon2008chaotic,lundstrom2008fractional,werner2010fractals,thurner2003scaling,teich1997fractal}. Fractional-order systems have also been explored in the context of neurophysiological networks constructed from electroencephalographic (EEG), electrocorticographic (ECoG), or blood-oxygen-level-dependent (BOLD) data~\citep{chatterjee2020fractional,magin2006fractional}.

Furthermore, we provide an overview of the work that has been done on controlling, estimating, and predicting neural dynamical systems modeled using fractional-order dynamics both in the continuous-time and discrete-time domains, towards the next generation of CNS. Specifically, the focus of our brief survey is threefold:

\begin{itemize}
    \item \textbf{Control:} We review different methods of controlling fractional-order systems, including a few previously presented in~\citep{efe2011fractional}. The work in \citep{birs2019survey} presents a survey of recent advances in fractional order control for time delay systems, and the works in \citep{matuvsuu2011application,chen2009fractional} provide overviews of the application of fractional calculus to control theory. In this paper, we review proportional-integral control, sliding mode control, backstepping control, adaptive control, optimal control, and model predictive control for fractional-order systems. Control of fractional-order systems is important to study so as to develop methods and therapies to mitigate and potentially eliminate diseases in the brain. 
    \item \textbf{System Identification:} We focus on  how to estimate the parameters of fractional-order systems from brain measurements, such as electroencephalography (EEG) and electrocorticography (ECoG), which is a necessary step in understanding the intricacies of the brain.  Extensive attention has been paid to estimating the parameters of fractional-order systems, as evidenced in the following comprehensive works:~\citep{ljung1999sysid} and~\citep{soderstrom1988system}. 
    \item \textbf{Estimation:}  We discuss the methods for estimating and predicting the state of fractional-order systems~\citep{sabatier2012observability,sierociuk2006fractional,Safari1,Safari2,miljkovic2017ecg,najar2009discrete,chatterjee2019dealing}. This problem is important in anticipating and mitigating irregular brain behavior such as an oncoming seizure. \citep{chatterjee2021minimum} proposed the design of a minimum-energy estimation framework for discrete-time fractional-order networks, where they assume that the state and output equations are affected by an additive disturbance and noise, respectively, that is considered to be deterministic, bounded, and unknown. First proposed by Mortensen~\citep{mortensen1968maximum}, and later refined by Hijab~\citep{hijab1979minimum}, minimum-energy estimators produce an estimate of the system state that is the ``most consistent" with the dynamics and the measurement updates of the system~\citep{fleming1997deterministic,willems2002deterministic,buchstaller2020deterministic,swerling1971modern,bonnabel2014contraction,fagnani1997deterministic,krener2003convergence,aguiar2006minimum,hassani2009multiple,pequito2009entropy,alessandretti2011minimum,ha2018cooperative,haring2020stability}. 
\end{itemize}

In what follows, we overview some of the definitions of both continuous and discrete fractional-order dynamics, towards assessing and designing different aspects of modeling, analysis, and control of the aforementioned. We start by introducing the fractional-order system in continuous-time \citep{kilbas2001differential,kilbas2002differential}.  

\subsection{Continuous-time fractional order systems}\label{sec:control_of_cont_frac} 

Riemann-Liouville and Caputo proposed the two popular (and equivalent) definitions of fractional order differintegration. 
Caputo's definition is the one widely used in control systems engineering, and it is the following: 
\begin{equation}\label{eq:caputo}
\begin{aligned}
\Delta^{\alpha} \sigma(t) &=\frac{1}{\Gamma(m-\alpha)} \int_{0}^{t} \frac{\Delta^{m} \sigma(\tau)}{(t-\tau)^{\alpha+1-m}} \mathrm{~d} \tau, 
\end{aligned}
\end{equation}
where $\alpha \in \mathbb{R}^{+}$is the differentiation order \citep{baleanu2012fractional}. 
Given the definition in~\eqref{eq:caputo}, let $m\in\mathbb{Z}$ with $m-1<\alpha<m$. For an $m$ satisfying the previous relation, the $\alpha$ order derivative of a function of time, $\sigma(t)$, has the following Laplace transform:

\begin{equation}\label{eq:caputo_2}
\begin{aligned}
\int_{0}^{\infty} e^{-s t} \Delta^{\alpha} \sigma(t) \mathrm{~d} t &=s^{\alpha} S(s)-\sum_{k=0}^{m-1} s^{\alpha-k-1} \mathrm{D}^{k} \sigma(0),
\end{aligned}
\end{equation}
where $\Gamma(\alpha)=\int_{0}^{t} e^{-t} t^{\alpha-1} \mathrm{~d} t$ is the Gamma function and $S(s)=\int_{0}^{\infty} e^{-s t} \sigma(t) \mathrm{~d} t$.

If a system is, initially, in a resting state (i.e., the initial conditions are zero), then the operator $\Delta^{\alpha}$ acting in the time domain has a counterpart $s^{\alpha}$ in the $s$-domain. 
In this case, we can describe the transfer function of a system by a fractional order differential equation 
\begin{equation}\label{eq:frac_init_rest}
    \begin{aligned}
&\left(a_{n} \Delta^{\alpha_{n}}+a_{n-1} \Delta^{\alpha_{n-1}}+\ldots+a_{1} \Delta^{\alpha_{1}}+a_{0}\right) y(t) \\
&\quad=\left(b_{m} \Delta^{\beta_{m}}+b_{m-1} \Delta^{\beta_{m-1}}+\ldots+b_{1} \Delta^{\beta_{1}}+b_{0}\right) u(t),
\end{aligned}
\end{equation}
which we can obtain as 
\begin{equation}\label{eq:frac_init_rest_2}
    \begin{aligned}
\frac{Y(s)}{U(s)}=\frac{a_{n} s^{\alpha_{n}}+a_{n-1} s^{\alpha_{n-1}}+\ldots+a_{1} s^{\alpha_{1}}+a_{0}}{b_{m} s^{\beta_{m}}+b_{m-1} s^{\beta_{m-1}}+\ldots+b_{1} s^{\beta_{1}}+b_{0}},
\end{aligned}
\end{equation}
where $a_{k}, b_{k} \in \mathbb{R}$ and $\alpha_{k}, \beta_{k} \in \mathbb{R}^{+}$.

Therefore, we can write an affine and fractional-order nonlinear system~as
\begin{equation}\label{eq:nonlin_frac_ord}
\Delta^{\alpha} x(t)=\mathbf{f}(x(t))+\mathbf{g}(x(t)) u(t)
\end{equation}
where $u$ is the control input, and $\mathbf{f}$ and $\mathbf{g} \neq 0$ are the vector functions of the system-state $x(t)$. 
If the system is linear, then we can describe its state-space representation as
\begin{equation}\label{eq:lin_frac_ord}
\begin{aligned}
\Delta^{\alpha} x(t) &=A x(t)+B u(t), \quad x(t) \in \mathbb{R}^{n} \\
y(t) &=C x(t)+D u(t).
\end{aligned}
\end{equation}
The transfer function characterizing the relation between $Y(s)$ and $U(s)$, and the Laplace transforms of the output and input, respectively, is as given as
\begin{equation}\label{eq:transf_fun}
    H(s) =C\left(s^{\alpha} I-A\right)^{-1} B+D.
\end{equation}
Therefore, we have the following solution for the homogeneous case $(u(t)=0)$
\begin{equation}
    x(t) =E_{\alpha}\left(A t^{\alpha}\right) x(0)=\Phi(t) x(0),
\end{equation}
where $E_{\alpha}\left(A t^{\alpha}\right)$ is the \emph{Mittag-Leffler function}~\citep{oldham1974fractional} defined as 
\[
\Phi(t)\equiv E_{\alpha}
\left(A t^{\alpha}\right)=\sum_{k=0}^{\infty}
\left(\frac{\left(A t^{\alpha}\right)^{k}}{\Gamma(1+\alpha k)}\right). 
\] 

Therefore, we can write the solution of the fractional state equation and the output equation in~\eqref{eq:lin_frac_ord} as
\[
y(t)=C \Phi\left(t-t_{0}\right) x\left(t_{0}\right)+C \displaystyle\int_{0}^{t} \Phi(t-\tau) B u(\tau) \mathrm{~d} \tau+D u(t).
\]

Despite the fact that real-world systems have continuous-time signals in nature, in reality, we measure and control these systems using digitized technologies, which motivates the study of discrete-time fractional-order systems \citep{caponetto2010fractional,goodrich2015discrete,mahmoud2012advances}. Subsequently, we now introduce the discrete-time description of the fractional-order dynamics.

\subsection{Discrete-time fractional-order systems}

In what follows next, we briefly introduce linear discrete-time \mbox{fractional-order} system models. A discrete-time linear fractional-order system model is described as follows:
\begin{equation}
\label{eq:fosmain}
    \Delta^\alpha x[k+1] = Ax[k] + Bu[k] + w[k],
\end{equation}
where $x[k]$ is the state for time step $k \in \mathbb{N}$, $A\in\mathbb{R}^{n\times n}$ is the state coupling matrix and \mbox{$\alpha \in (\mathbb{R}^n)^{+}$} is the vector of fractional-order coefficients. The signal $u[k] \in \mathbb{R}^{n_u}$ denotes the input corresponding to the actuation signal and the matrix $B \in \mathbb{R}^{n \times n_u}$ is the input matrix that scales the actuation signal. The term $w[k]$ denotes the process noise or additive disturbance, whose stochastic characterization (or the lack thereof) will be clear from the context in which these systems are being used. These models are similar to classical \mbox{discrete-time} linear \mbox{time-invariant} system models with the exception of the inclusion of the (Gr\"unwald-Letnikov) fractional derivative, whose expansion and discretization for the $i$-th state, $1 \leq i \leq n$, can be expressed as
\begin{equation}
\label{eq:GL_deriv}
    \Delta^{\alpha_i} x_i[k] = \sum_{j=0}^k \psi (\alpha_i,j) x_i[k-j],
\end{equation}
where $\alpha_i$ is the fractional-order coefficient corresponding to the state $i$ and
\begin{equation}
\label{eq:weight_frac}
    \psi(\alpha_i,j) = \frac{\Gamma(j-\alpha_i)}{\Gamma(-\alpha_i) \Gamma(j+1)},
\end{equation}
with $\Gamma(\cdot)$ being the gamma function defined by $\Gamma(z) = \int_0^{\infty} s^{z-1} e^{-s} \: \mathrm{d}s$ for all complex numbers $z$ with positive real part, $\mathfrak{R}(z) > 0$~\citep{vinagre2000some,DzielinskiFOS}. Simply put, larger values of the \mbox{fractional-order} coefficients imply a lower dependency on the previous data from that state (i.e., a faster decay of the weights used as linear combination of previous data).

We now review some essential theory for \mbox{fractional-order} systems, including an approximation of~\eqref{eq:fosmain} with $u[k]=0$ for all $k \in \mathbb{N}$ as an LTI system. Using the expansion of the Gr\"unwald-Letnikov derivative in~\eqref{eq:GL_deriv}, we have
\begin{align}
    \Delta^\alpha x[k] &= \begin{bmatrix}\Delta^{\alpha_1}x_1[k]\\ \vdots \\ \Delta^{\alpha_n}x_n[k] \end{bmatrix} = \begin{bmatrix}\sum_{j=0}^k\psi(\alpha_1,j)x_1[k-j]\\ \vdots \\ \sum_{j=0}^k\psi(\alpha_n,j)x_n[k-j]\end{bmatrix} \nonumber \\ &= \sum_{j=0}^k\underbrace{\begin{bmatrix}\alpha_1 & \ldots & 0\\ \vdots & \ddots & \vdots \\ 0 & \ldots & \alpha_n\end{bmatrix}}_{D(\alpha,j)}\begin{bmatrix}x_1[k-j] \\ \vdots\\ x_n[k-j]\end{bmatrix} \nonumber \\ &= \sum_{j=0}^k D(\alpha,j)x[k-j].
    \label{eq:Deltaalpha}
\end{align}
The above formulation distinctly highlights one of the main peculiarities of DT-FODS in that the fractional derivative $\Delta^{\alpha} x[k]$ is a weighted linear combination of not just the previous state but of every single state up to the current one, with the weights given by~\eqref{eq:weight_frac} following a power-law decay.

Plugging~\eqref{eq:Deltaalpha} into the DT-FODS formulation~\eqref{eq:fosmain} with $u[k]=0$ for all $k \in \mathbb{N}$, we have
\begin{equation}
    \sum_{j=0}^{k+1} D(\alpha,j) x[k+1-j] = Ax[k] + w[k],
\end{equation}
or, equivalently,
\begin{equation}
    D(\alpha,0) x[k+1] = - \sum_{j=1}^{k+1} D(\alpha,j) x[k+1-j] + Ax[k] + w[k],
\end{equation}
which leads to
\begin{equation}
\label{eq:D_evol}
    x[k+1] = -\sum_{j=0}^k D(\alpha,j+1) x[k-j] + Ax[k] + w[k],
\end{equation}
since $D(\alpha,0) = I_n$, where $I_n$ is the $n \times n$ identity matrix. Alternatively,~\eqref{eq:D_evol} can be written as
\begin{align}\label{eq:state_evol_2}
    x[k+1] &= \sum_{j=0}^k A_j x[k-j] + w[k] \nonumber \\
    x[0] &= x_0,
\end{align}
where
\begin{equation}
    A_j = \begin{cases}A -\textnormal{diag}(\alpha_1,\ldots,\alpha_n) & \text{if} \; j=0\\ -D(\alpha,j+1) & \text{if} \; j \geq 1 \end{cases}.
\end{equation}

\subsection{Paper organization}

Section \ref{sec:stability} presents the results on the stability of fractional-order systems. Section \ref{sec:contr_obsv} provides a summary of the work on controllability and observability of fractional-order systems. Section~\ref{sec:pid_control} summarizes the work on \mbox{proportional-integral-derivative} controllers for fractional-order dynamical systems. Section~\ref{sec:sliding} reviews sliding mode control for fractional-order systems. Section~\ref{sec:backstep} outlines the procedure for constructing a backstepping controller for fractional-order systems.  Section~\ref{sec:adaptive} summarizes adaptive control for fractional-order systems. Section~\ref{sec:system_identification} discusses the methods for performing system identification on fractional-order systems. Section~\ref{sec:state_estimation} overviews the techniques for state estimation of fractional-order systems, including the method known as minimum-energy state estimation. Section \ref{sec:fractional_optimal_control} presents fractional optimal control for continuous-time fractional-order systems. Section~\ref{sec:model_predictive_control} gives a background on model predictive control for \mbox{fractional-order} systems. Section~\ref{sec:applications} summarizes simulation results pertaining to system identification, state estimation, and closed-loop control of fractional cyber-neural systems. Finally, section~\ref{sec:future_research} presents possible directions for future research.

\section{Stability}\label{sec:stability}
Stability can be described as the behavior of the state of a system after a reasonable amount of time. While there are different notions of stability, in affect, a system is stable if the behavior of the system is bounded. The prior literature describes conditions for continuous-time fractional-order systems~\citep{benzaouia2014stabilization, li2009mittag,monje2010fractional} and for single-input single-output continuous-time  \emph{commensurate} systems (i.e., systems with equal fractional exponents across state variables)~\citep{dastjerdi2019linear}. \citep{li2010stability} provides the generalized Mittag–Leffler stability conditions of continuous-time fractional-order systems using the Lyapunov direct method. In what follows, we summarize the stability conditions for continuous-time commensurate fractional-order systems. 
Let $\sigma(A)=\{\lambda_1,\ldots,\lambda_n\}$ be the spectrum (set of eigenvalues) of $A$. 
We say that the commensurate system in~\eqref{eq:lin_frac_ord} is stable if 
\begin{equation}\label{eq:stab}
    \arg \left(\lambda_{i}\right)>\alpha \frac{\pi}{2}, \quad\text{for all } i=1,\ldots, n,
\end{equation}
where $\arg(z)$, in the complex plane, is the 2D polar angle $\varphi$
 from the positive real axis to the vector representing $z$, and $0<\alpha<2$ \citep{rivero2013stability}. 
In the case of the transfer function in~\eqref{eq:transf_fun}, we have that $\sigma(A)$ corresponds to the poles of the system, and the previous stability condition of~\eqref{eq:stab} also applies. 
Notice that, in the integer order case $(\alpha=1)$, the stability condition of~\eqref{eq:stab} describes the open left half $s$-plane. 
For a more detailed discussion on the stability of continuous-time systems, we refer the reader to~\citep{matignon1996stability,chen2005robust,ortigueira2000introduction}.

For discrete-time fractional-order systems, the authors of \citep{dzielinski2008stability} leverage an infinite dimensional representation of truncated discrete-time fractional-order systems (i.e., with finite memory) to give conservative sufficient conditions for stability. While the work in \citep{buslowicz2013necessary} does provide necessary and sufficient conditions for practical and asymptotic stability of discrete-time fractional-order systems, they only consider commensurate-order systems (i.e., $\alpha$ is the same for all state variables). Recent work has introduced stability conditions for multivariate discrete-time fractional-order systems with arbitrary fractional exponents and leverages these conditions to study the stability of a real-world EEG cognitive motor data set modeled as a discrete-time fractional-order system and to provide evidence of its relevance in the context of cognitive motor control \citep{Reed2021}. 

That said, a simple to state necessary and sufficient condition like~\eqref{eq:stab} for non-commensurate both continuous and discrete-time systems is still missing. This limits the capacity to assess the stability of such systems and their applicability in the context of neural systems and possibly some neurological diseases such as epilepsy.

\section{Controllability and observability} \label{sec:contr_obsv}
Controllability is a prerequisite in the ability to manipulate a system to any desired state in a finite amount of time. On the other hand, observability is necessary to obtain a complete picture of the system on the whole. For continuous-time systems, \citep{matignon1996some} gives results on the controllability and observability of finite-dimensional continuous-time fractional-order systems. \citep{balachandran2013observability} gives a comprehensive overview of the conditions for controllability and observability of continuous-time linear fractional-order systems.  Similarly, \citep{guermah2008controllability} provides these results for discrete-time linear fractional-order systems. The work in  \citep{mozyrska2012fractional} derives the conditions for controllability and observability of finite memory discrete-time fractional-order systems.

Previous work has examined the design of controllable networks exhibiting discrete-time linear fractional-order dynamics using energy-based methods  \citep{kyriakis2020effects} and by maximizing the rank of the controllability matrix through a greedy algorithm \citep{cao2019actuation}.  Similarly, there has been work in selecting the minimal number of EEG sensors to achieve observability for discrete-time fractional-order systems \citep{gupta2018dealing,Xue2016,tzoumas2018selecting,PequitoC24}.    

A system is \emph{controllable} if there exists a control input such that the final state can be driven to zero in a finite amount of time. In particular, for continuous-time linear fractional-order systems modeled by \eqref{eq:lin_frac_ord}, the system is controllable on $[t_0, t_1]$ if for every pair of vectors $x(t_0), x(t_1)\in\mathbb{R}^{n}$, there is a control $u(t)\in L^{2}([t_0,t_1],\mathbb{R}^{m})$ such that the solution $x(t)$ of \eqref{eq:lin_frac_ord} which satisfies $x(t_0)=x_0$ also satisfies $x(t_1)=x_1$, where $L^{2}([t_0,t_1],\mathbb{R}^{m})$ is the space of all square integrable $\mathbb{R}^{m}$ valued measureable functions defined on $[t_0, t_1]$. Thus, we say that \eqref{eq:lin_frac_ord} is controllable on $[t_0,t_1]$ if and only if the controllability Gramian matrix $$ \int_{t_0}^{t_1}(t_1-\tau)^{\alpha -1}E_{\alpha,\alpha}(A(t_1 - \tau)^{\alpha})BB^{\intercal}E_{\alpha,\alpha}(A^{\intercal}(t_1-\tau)^{\alpha})\mathrm{~d}\tau$$ is positive definite for some $t_1>t_0$ (Theorem 3, \citep{balachandran2013observability}).

For discrete-time linear fractional-order
system modeled by \eqref{eq:fosmain}, the system is controllable if there exists a control sequence $\{\mathbf{u}[0],\ldots,\mathbf{u}[T-1]\}$ such that $\mathbf{x}[T] = \mathbf{0}$
from any initial state $\mathbf{x}[0]\in\mathbb{R}^{n}$ in a finite time \citep{guermah2008controllability}. To present the conditions for controllability and observability for discrete-time fractional-order systems, we first start by noticing that the discrete-time linear fractional-order system~\eqref{eq:fosmain} can be re-written as~\citep[Lemma 2]{gupta2018dealing}:
\begin{equation} 
x[k] = G_{k}x[0],
\label{eq:G}
\end{equation}
where 
\begin{align}\label{eq:state_transition_matrix_G}
G_{k} = 
\begin{cases}
I_{n}, & k=0 \\
\sum_{j=0}^{k-1}A_{j}G_{k-1-j}, & k \geq 1
\end{cases}
\end{align}
with $A_{0} = A-D(\alpha,1)$, $A_{j} = -D(\alpha, j+1)$, for  $j \geq 1$, and 
\begin{equation}
D(\alpha,j) = \begin{bmatrix} \psi(\alpha_{1},j) & 0 & \dots & 0 \\
0 & \psi(\alpha_{2},j) & \dots & 0 \\ 
0 & \vdots & \ddots & 0 \\
0 & 0 & \dots & \psi(\alpha_{n},j)
\end{bmatrix}.
\end{equation}

The linear discrete-time fractional-order
system modeled by \eqref{eq:fosmain} is controllable if and only if there exists a finite time $K$ such that $\text{rank}(W_c(0,K)) =
n$, where $W_c(0,K) = G_{K}^{-1}\sum_{j=0}^{K-1}G_{j}BB^{\intercal}G_{j}^{\intercal}G_{K}^{-\intercal}$~\citep[Theorem 4]{guermah2008controllability}. Furthermore, an input sequence $\begin{bmatrix}\mathbf{u}^{\intercal}[K-1], \mathbf{u}^{\intercal}[K-2], \dots \mathbf{u}^{\intercal}[0]\end{bmatrix}^{\intercal}$ that transfers $\mathbf{x}[0] \neq 0$ to $\mathbf{x}[K] = 0$ is given by
\begin{equation}
   \begin{bmatrix}\mathbf{u}[K-1]\\
    \mathbf{u}[K-2]\\  
    \vdots \\ 
    \mathbf{u}[0]\end{bmatrix}= -[G_0B G_1B \ldots G_{K-1}B]^{\intercal}G_{K}^{-\intercal} W_c^{-1}(0,K)\mathbf{x}[0].
\end{equation}

Similarly, a system is \emph{observable} if and only if the initial state $\mathbf{x}[0]$ can be uniquely determined from the knowledge of the control input and observations. For continuous-time systems, the system is observable on an interval $[t_0,t_1]$ if $y(t)=Cx(t)=0$ for $t\in[t_0,t_1]$ implies $x(t)=0$ for $[t_0,t_1].$ Hence, the system in \eqref{eq:lin_frac_ord} is observable on $[t_0,t_1]$ if an only if the observability Gramian matrix $W=\int_{t_0}^{t_1}E_{\alpha}(A^{\intercal}(t-t_0)^{\alpha})C^{\intercal}CE_{\alpha}(A(t-t_0)^{\alpha})\mathrm{~d}t$ is positive definite~\citep[Theorem 1]{balachandran2013observability}. 

For linear discrete-time fractional-order systems modeled by \eqref{eq:fosmain}, the system is said to be observable if and only if there exists some $K>0$ such that the initial state $\mathbf{x}[0]$ at time $k=0$ can be uniquely determined from the knowledge of $\{\mathbf{u}[0],\ldots,\mathbf{u}[K-1]\}$ and $\{\mathbf{y}[0],\dots,\mathbf{y}[K-1]\}$. Therefore, by Theorem~5 in \citep{guermah2008controllability}, the linear discrete-time fractional-order system is observable if and only if there exists a finite time $K$ such that $\text{rank}(\mathcal{O}_K)=n$, where $\mathcal{O}_K=\begin{bmatrix} CG_0, CG_1, \ldots, CG_{K-1} \end{bmatrix}^{\intercal}$ or, equivalently, $\text{rank}(W_o(0,K))=n$, where $W_o(0,K)=\sum_{j=0}^{K-1}G_{j}^{\intercal}C^{\intercal}CG_{j}$. Furthermore, the initial state at $\mathbf{x}[0]$ is given by 
\begin{equation}
    \mathbf{x}[0] = W_{o}^{-1}(0,K)\mathcal{O}_{K}^{\intercal}[\tilde{\mathcal{Y}}_K-\mathcal{M}_K\tilde{\mathcal{U}}_K],
\end{equation}

\noindent where $\tilde{\mathcal{U}}_K = \begin{bmatrix} \mathbf{u}^\intercal[0], \mathbf{u}^{\intercal}[1], \dots, \mathbf{u}^{\intercal}[K-1] \end{bmatrix}$, $\tilde{\mathcal{Y}}_K =\begin{bmatrix} \mathbf{y}^{\intercal}[0], \dots,  \mathbf{y}^{\intercal}[K-1] \end{bmatrix}^\intercal$, and 
\begin{equation*}
    \mathcal{M}_K = \begin{bmatrix} 
0 & 0 & \dots & 0 & 0 \\
CG_{0}B & 0 & \dots & 0 & 0 \\
CG_{1}B & CG_{0}B & \dots & 0 & 0 \\
CG_{2}B & CG_{1}B & \dots & 0 & 0 \\
\vdots & \vdots & \ddots & \vdots & \vdots \\ 
CG_{K-2}B & CG_{K-3}B & \dots & CG_{0}B & 0
\end{bmatrix}.
\end{equation*}

\section{Proportional-integral-derivative control}\label{sec:pid_control}


A \emph{proportional-integral-derivative controller} (\emph{PID controller} or \emph{three-term controller}) is a control loop mechanism. 
It employs feedback and is commonly utilized in industrial control systems and a variety of other applications that need a continuously modulated control. 
A PID controller continuously computes an error value as the difference between the desired setpoint (SP) and a measured process variable (PV) and implements a correction based on proportional, integral, and derivative terms (P, I, and D, respectively).

A simple and practical example is the cruise control of a vehicle. 
For example, if only constant engine power is applied, a car ascending a hill would lose speed. 
In a situation like this, the controller's PID algorithm is responsible for restoring the measured speed to the desired cruise control speed. Moreover, the PID control increases the power output of the engine in a controlled manner, with minimal delay and overshoot. 

The interpretability and comprehensibility of PID controllers make it a typical choice regarding peripherals that automatically tune the system's parameters without external intervention. The work in \citep{podlubny1999fractional} provides the framework for fractional-order PID control, which we summarize next. 
We can define the fractional order version of PID controllers using the following transfer function:
\begin{equation}\label{eq:pid}
C(s)=k_{p}+\frac{k_{i}}{s^{\lambda}}+k_{d} s^{\mu}.
\end{equation}
For $\lambda=1$ and $\mu=1$, we obtain the standard integer order setting with three degrees of freedom: $k_{p}, k_{i}$, and $k_{d}$. 
Notwithstanding, in~\eqref{eq:pid}, we have five parameters to determine, yielding five independent specifications that we can force to meet. 
If we place the controller in front of a $G(s)$ in a unity feedback loop, then the first specification can be on phase margin as it is tightly coupled with the robustness of the control system. 
The equations that define the phase margin are $20 \log \left|C\left(w_{q c}\right) G\left(w_{q c}\right)\right|=0 \mathrm{~dB}$ and $\arg \left(C\left(w_{g c}\right) G\left(w_{g c}\right)\right)=-\pi+\varphi_{p m}$, where $w_{g c}$ is the gain crossover frequency and $\varphi_{p m}$ the phase margin. 

Subsequently, we may force a flat magnitude response $|G(j w) C(j w)|$ around the gain crossover frequency. 
We can ensure this response, see~\citep{monje2008tuning}, by setting the derivative $\frac{\mathrm{d}}{\mathrm{d}w}(\arg (C(j w) G(j w)))$ to zero when $w=w_{c g}$. 
Moreover, ensuring this constraint makes the closed-loop control system robust against variations in the gain of $G(s)$. 

Another specification supposes that a controller introduces the property of noise rejection in high frequencies.  
This property can be accomplished by fixing a critical frequency, $w_{h}$.  If this frequency is exceeded, then the magnitude of the transfer function $T=C G /(1+C G)$ (corresponding to the complementary sensitivity function) is smaller than a preselected level. 

Next, how can we ensure that a good output disturbance is not rejected? 
To address this, we can force an upper bound $(M)$ on the sensitivity function's magnitude below a predefined frequency $\left(w_{s}\right)$. 
Hence, it follows that
\[
20 \log |S(j w)|_{w \leq w_{s}}=20 \log \left|\frac{1}{1+C(j w) G(j w)}\right|_{w \leq w_{s}}\leq M \mathrm{~dB}.
\]
In the final step, to achieve a zero steady-state error, we need to design the controller $C(s)$ with an integral component.

Notice that, although solving the necessary set of equations from the constraints above is one way to establish the parameters $k_{p}, k_{i}, k_{d}, \lambda$, and $\mu $, this requires the prior knowledge of model order, dead time, poles and zeros. 
In the scenario where we do not have this prior knowledge, we may, alternatively, resort to \emph{autotuning}~\citep{monje2008tuning,chen2004robust}.

\section{Sliding mode control}\label{sec:sliding}

In control systems, \emph{sliding mode control} (SMC) is a nonlinear control method that adjusts the dynamics of a nonlinear system by applying a discontinuous control signal (a set-valued control signal). 
This control signal compels the system to ``slide'' along a cross-section of the system's normal behavior. 
In this case, the state-feedback control law is not a function continuous in time. 
Instead, it can switch between continuous structures based on the state space current position to achieve the desired behavior. 
There are two stages in SMC:  
\emph{(i)} the \emph{reaching phase}, which is the phase that lasts until the hitting of a trajectory to the switching subspace;
\emph{(ii)} the \emph{sliding mode}, which is the motion after the previous phase. 
A relevant property of stage \emph{(ii)} is the robustness against disturbances and variations in the process parameters -- i.e., the \emph{invariance property}. 

Now, we present a set of results regarding SMC for fractional order systems. 
Given the $n$th order fractional dynamic system in~\eqref{eq:nonlin_frac_ord} and the following switching function
\begin{equation}\label{eq:Lambda}
\sigma(t)=\Lambda\left(c(t)-r(t)\right),
\end{equation}
where $\Lambda$ is a parameter designed to make the sliding manifold defined by $\sigma=0$ to be a stable subspace, where the stability can be settled via~\eqref{eq:stab}. 
This entails that, despite the process being nonlinear, the nominal plant model is linear. 
If $0<\alpha<1$ and $r$ is the vector of differentiable command signals, then the goal of the reaching law approach is to get $\Delta^{\alpha} \sigma(t)=-k \operatorname{sgn}(\sigma(t))$ for some $k>0$. 
When $\alpha=1$, it corresponds to $\dot{\sigma}(t)=-k \operatorname{sgn}(\sigma(t))$, which ensures $\sigma(t) \dot{\sigma}(t)<0$ whenever $\sigma \neq 0$.
This solution is the time derivative of the Lyapunov function $V=\frac{1}{2} \sigma(t)^{2}$, which physical meaning is to provide the sliding manifold an attractor such that, once the error vector gets trapped to it, the subsequent motion takes place in the proximity of the sliding hypersurface.

Next, we need to show that the aforementioned mechanism also works for non-integer differentiation order~\citep{vinagre2006fractional}. 
We start by differentiating $\Delta^{\alpha} \sigma(t)=-k \operatorname{sgn}(\sigma(t))$ at the order $-\alpha$
\[
\Delta^{1}\left(\Delta^{-\alpha}\left(\Delta^{\alpha} \sigma(t)\right)\right) =-k \Delta^{1}\left(\Delta^{-\alpha} \operatorname{sgn}(\sigma(t))\right) 
\]
and, next, differentiate at order unity to obtain $\dot{\sigma}(t)$ 
\[
\dot{\sigma}(t) =-k \Delta^{1-\alpha} \operatorname{sgn}(\sigma(t)).
\]
Because $0<\alpha<1$, it follows that $\operatorname{sgn}\left(\Delta^{1-\alpha} \operatorname{sgn}(\sigma(t))\right)=\operatorname{sgn}(\sigma(t))$. 
Forcing $\Delta^{\alpha} \sigma(t)=-k \operatorname{sgn}(\sigma(t))$ makes the locus described by $\sigma=0$ a global attractor. 

It is easy to check that choosing $\Delta^{\alpha} \sigma(t)=-k \operatorname{sgn}(\sigma(t))-p \sigma(t)$ with $p>0$ has the same effect on the reaching dynamics of that in the integer order design. 
Notice that with $p \sigma=p|\sigma(t)| \operatorname{sgn}(\sigma)$, the following relation holds between $\dot{\sigma}(t)$ and $\operatorname{sgn}(\sigma(t))$:
\[
\begin{aligned}
\dot{\sigma}(t) &=-k \Delta^{1-\alpha} \operatorname{sgn}(\sigma(t))-p \Delta^{1-\alpha}(|\sigma(t)| \operatorname{sgn}(\sigma(t))) \\
&=-\Delta^{1-\alpha}((k+p|\sigma(t)|) \operatorname{sgn}(\sigma(t))).
\end{aligned}
\]
Notice that, since $\operatorname{sgn}\left(\Delta^{1-\alpha} \operatorname{sgn}(\sigma(t))\right)=\operatorname{sgn}(\sigma)$, the reaching dynamics governed by the above expression generates a stronger push from both sides of the switching manifold. This effect translates into the attraction strength of the switching manifold being higher, for any $\sigma(t)$ with $p \neq 0$, than for $p=0$.
Moreover, for a fixed $\sigma(t)$, larger values of $p$ create larger values of $\dot{\sigma}(t)$, which leads to reaching quicker the place characterized by $\sigma=0$.  
If we select the Lyapunov function $V(t)=\frac{1}{2} \sigma(t)^{2}$ and compute its $\alpha$th order derivative, using the Leibniz's differentiation rule, we obtain  
$$\Delta^{\alpha} V(t)=\sum_{k=0}^{\infty}\frac{\Gamma(1+\alpha)}{\Gamma(1+k)} \Gamma(1-k+\alpha) \Delta^{k} \sigma(t) \Delta^{\alpha-k} \sigma(t),$$ 
i.e., an expression with infinitely many terms. 
Therefore, we are not able to infer the attractiveness of $\sigma(t)=0$, deduced from $\sigma(t) \Delta^{\alpha} \sigma(t)<0$, or more specifically, from  $\Delta^{\alpha} \sigma(t)=-k \operatorname{sgn}(\sigma(t))-p \sigma(t)$. 

Recalling definition~\eqref{eq:caputo}, the following equality holds  
$$
\sigma(t) \Delta^{\alpha} \sigma(t)=\frac{\sigma(t)}{\Gamma(1-\alpha)} \int_{0}^{t}\frac{\Delta \sigma(\tau) }{(t-\tau)^{\alpha}} \mathrm{d} \tau.
$$ 
The previous relation imposes two possibilities to have that $\sigma(t) \Delta^{\alpha} \sigma(t)<0$: 
\begin{itemize}
\item [(i)]if $\sigma(t)>0$, then $\Delta \sigma(t)$ (the first derivative of $\sigma(t)$) must be negative;
\item[(ii)] if $\sigma(t)<0$, then $\Delta\sigma(t)$ (the first derivative of $\sigma(t)$) must be positive.
\end{itemize} 
In conclusion, an appropriately designed control law is sufficient for closed-loop stability, forcing $\sigma \Delta^{\alpha} \sigma(t)<0$. 
Therefore, the stability requirement $\sigma(t) \dot{\sigma}(t)<0$ (or $\sigma(t) \Delta \sigma(t)<0$) of the integer order design is obtained naturally, whenever we impose $\sigma(t) \Delta^{\alpha} \sigma(t)<0$. 

In~\citep{efe2011fractional}, the following is proposed. 
Compute $\alpha$th order derivative of~\eqref{eq:Lambda}, which is
$$
\Delta^{\alpha} \sigma(t)=\mathbf{\Lambda}\left(\Delta^{\alpha} x(t)-\Delta^{\alpha} r(t)\right)=\mathbf{\Lambda}\left(\mathbf{f}(x(t))+\mathbf{g}(x(t)) u(t)-\Delta^{\alpha} r(t)\right).
$$
Next, setting the previous expression equal to $-k \operatorname{sgn}(\sigma(t))-p \sigma(t)$ and solving for $u$ yields the following control signal:

\begin{equation}\label{eq:control_sig}
u(t)=\frac{-\Lambda \mathbf{f}(x(t))+\Lambda \Delta^{\alpha} r(t)-k \operatorname{sgn}(\sigma(t))-p \sigma(t)}{\Lambda g(x(t))},
\end{equation}
where we need to have that $\Lambda g(x(t)) \neq 0$. 
Having the encountered control law, deduced from a nominal model, an important question is what would be the response of the system, whenever the model in~\eqref{eq:nonlin_frac_ord} is a nominal representation of a plant holding the uncertainties $\Delta \mathbf{f}(x(t))$ and $\Delta \mathbf{g}(x(t))$, such as the following
\begin{equation}\label{eq:uncert_eq}
\Delta^{\alpha} x(t)=(\mathbf{f}(x(t))+\Delta \mathbf{f}(x(t)))+(\mathbf{g}(x(t))+\Delta \mathbf{g}(x(t))) u(t).
\end{equation}
Combining~\eqref{eq:control_sig} and~\eqref{eq:uncert_eq} yields the following dynamics
\begin{equation}\label{eq:24}
\begin{split}
\Delta^{\alpha} \sigma(t)
=
-\left(1+\frac{\mathbf{\Lambda} \Delta\mathbf{g}(x(t))}{\boldsymbol{\Lambda} \mathbf{g}(x(t))}\right)(k \operatorname{sgn}(\sigma(t))+p \sigma(t)) \\
+\frac{\boldsymbol{\Lambda} \Delta \mathbf{g}(x(t))}{\boldsymbol{\Lambda} \mathbf{g}(x(t))} \boldsymbol{\Lambda}\left(\Delta^{\alpha} r(t)-\mathbf{f}(x(t))\right)+\boldsymbol{\Lambda} \Delta \mathbf{f}(x(t)).
\end{split}
\end{equation}
Hence, we have the following properties:
\begin{itemize}
    \item If there are no uncertainties, i.e., $\Delta \mathbf{f}(x(t))=\Delta \mathbf{g}(x(t))=$ 0, then we have $\Delta^{\alpha} \sigma(t)=-k \operatorname{sgn}(\sigma(t))-p \sigma(t)$, being desired to observe the sliding regime, after hitting the sliding hypersurface;
    \item If $\Delta \mathbf{g}(x(t))=0$ and the columns of $\Delta \mathbf{f}(x(t))$ are in the range space of $\mathbf{g}(x(t))$, then $\Delta^{\alpha} \sigma(t)=-k \operatorname{sgn}(\sigma(t))-p \sigma+\boldsymbol{\Lambda} \Delta \mathbf{f}(x(t))$.  This case further requires the hold of the condition in of $|A \Delta \mathbf{f}(x(t))|<k$ to ensure that $\sigma(t) \Delta^{\alpha} \sigma(t)<0$;
    \item If the uncertainty terms are nonzero, then~\eqref{eq:24} is valid, which implies that the designer has to carefully set $k$ and $p$ to keep the attractiveness of the subspace defined by $\sigma(t)=0$ The following two conditions are required to ensure that $\sigma(t) \Delta^{\alpha} \sigma(t)<0$:
\end{itemize}
    \[
\begin{array}{rcl}
\displaystyle\left|\frac{\Lambda \Delta g(x(t))}{\Lambda g(x(t))}\right| & < & 1 \\[.3cm]
\displaystyle\left(1+\frac{\Lambda \Delta g(x(t))}{\Lambda g(x(t))}\right) k & > & \displaystyle\left|\frac{\Lambda \Delta g(x(t))}{\Lambda g(x(t))} \Lambda\left(\Delta^{\alpha} r(t)-f(x(t))\right)+\Lambda \Delta f(x(t))\right|.
\end{array}
\]
The columns of $\Delta \mathbf{f}(x(t))$ and $\Delta \mathbf{g}(x(t))$ are assumed to be in the range space of $\mathbf{g}(x(t))$, i.e., the uncertainties are matched. 
If the previous condition is not satisfied, then the closed-loop performance will deteriorate.

Finally, notice that the first hitting to the switching subspace yields when $t=t_{h}$, where $t_{h}=(|\sigma(0)| \Gamma(\alpha+1) / k)^{1 / \alpha}$.

\section{Backstepping control}\label{sec:backstep}

\emph{Backstepping} is a technique developed in the 90s by Petar V. Kokotovic and others~\citep{kokotovic1992joy,lozano1992adaptive}. 
The goal of this technique is to design stabilizing controls for a special class of nonlinear dynamical systems. 
These systems consist of subsystems that radiate out from an irreducible subsystem, which we can stabilize using some method. 
Due to its recursive structure, the designer can start the design process at the \mbox{known-stable} system and ``back out'' new controllers that progressively stabilize each outer subsystem. 
The process of stabilization stops when the final external control is achieved. 
In other words, backstepping is based on the definition of a set of intermediate variables and the process of ensuring the negativity of Lyapunov functions that are combined to form a common control Lyapunov function for the overall system.


In fact, we can use the backstepping technique in a particular but wide class of systems. 
%
Consider the following system

\begin{equation}\label{eq:pos_vel}
\begin{aligned}
&x_{1}^{\left(\alpha_{1}\right)}(t)=x_{2}(t) \\
&x_{2}^{\left(\alpha_{2}\right)}(t)=\mathbf{f}\left(x_{1}(t), x_{2}(t)\right)+\mathbf{g}\left(x_{1}(t), x_{2}(t)\right) u(t),
\end{aligned}
\end{equation}
where $x_{1}$ and $x_{2}$ are the state variables, $0<\alpha_{1}, \alpha_{2}<1$ are positive fractional differentiation orders, $\mathbf{f}$ and $\mathbf{g}$ are known and smooth functions of the state variables such that $\mathbf{g}(x_1(t),x_2(t)) \neq 0$. 
Additionally, consider the intermediate variables of backstepping design:

\[
\begin{aligned}
&z_{1}(t)=x_{1}(t)-r_{1}(t)-A_{1}(t) \\
&z_{2}(t)=x_{2}(t)-r_{2}(t)-A_{2}(t),
\end{aligned}
\]
where $A_{1}(t)=0$ and $r_{1}^{\left(\alpha_{1}\right)}(t)=r_{2}(t)$.

Subsequently, consider the Lyapunov function with variable of interest $z$
\[
V(t)=\frac{1}{2} z^{2}(t).
\]
Now, from Section~\ref{sec:sliding}, we have that  $z(t) z^{(\alpha)}(t)$ ensures $z(t) \dot{z}(t)<0$, for any  $0<\alpha<1$. 
That said, we formulate the backstepping control technique for the plant described by~\eqref{eq:pos_vel}, by checking recurrently the quantities  
$z_{1}(t) z_{1}^{\left(\alpha_{1}\right)}(t)$ and $z_{1}(t) z_{1}^{\left(\alpha_{1}\right)}(t)+z_{2}(t) z_{2}^{\left(\alpha_{2}\right)}(t)$ as the following steps:
\begin{enumerate}
    \item  Check $z_{1}(t) z_{1}^{\left(\alpha_{1}\right)}(t)$:
\[
\begin{aligned}
z_{1}(t) z_{1}^{\left(\alpha_{1}\right)}(t) &=z_{1}(t)\left(x_{1}^{\left(\alpha_{1}\right)}(t)-r_{1}^{\left(\alpha_{1}\right)}(t)\right) \\
&=z_{1}(t)\left(x_{2}(t)-r_{2}(t)\right) \\
&=z_{1}(t)\left(z_{2}(t)+r_{2}(t)+A_{2}(t)-r_{2}(t)\right) \\
&=z_{1}(t)\left(z_{2}(t)+A_{2}(t)\right)
\end{aligned}
\]
\item Choose $A_{2}(t)=-k_{1} z_{1}(t)$, with $k_{1}>0$, this would entail that
\[
z_{1}(t) z_{1}^{\left(\beta_{1}\right)}(t)=-k_{1} z_{1}^{2}(t)+z_{1}(t) z_{2}(t)
\]
\item Check 
$
z_{1}(t) z_{1}^{\left(\alpha_{1}\right)}(t)+z_{2}(t) z_{2}^{\left(\alpha_{2}\right)}(t)
$:
\[
\begin{aligned}
&z_{1}(t) z_{1}^{\left(\alpha_{1}\right)}(t)+z_{2}(t) z_{2}^{\left(\alpha_{2}\right)}(t) \\
&\quad=-k_{1} z_{1}^{2}(t)+z_{1}(t) z_{2}(t)+z_{2}\left(x_{2}^{\left(\alpha_{2}\right)}(t)-r_{2}^{\left(\alpha_{2}\right)}(t)-A_{2}^{\left(\alpha_{2}\right)}(t)\right) \\
&\quad=-k_{1} z_{1}^2(t)+z_{2}(t)\left(x_{2}^{\left(\alpha_{2}\right)}(t)-r_{2}^{\left(\alpha_{2}\right)}(t)-A_{2}^{\left(\alpha_{2}\right)}(t)+z_{1}(t)\right)\\
&\quad=-k_{1} z_{1}^{2}(t)+z_{2}(t)\left(\mathbf{f}\left(x_1(t),x_2(t)\right)+\mathbf{g}\left(x_1(t),x_2(t)\right) u(t)\phantom{r_{2}^{\left(\alpha_{2}\right)}}\right.\\
&\qquad\qquad\qquad\qquad\qquad \left.-r_{2}^{\left(\alpha_{2}\right)}(t)-A_{2}^{\left(\alpha_{2}\right)}(t)+z_{1}(t)\right)
\end{aligned}
\]
\item Force $z_{1}(t) z_{1}^{\left(\alpha_{1}\right)}(t)+z_{2}(t) z_{2}^{\left(\alpha_{2}\right)}(t)=-k_{1} z_{1}^{2}(t)-k_{2} z_{2}^{2}(t)$, with $k_{2}>0$, which implies that
\[
\mathbf{f}(x_1(t),x_2(t))+\mathbf{g}(x_1(t),x_2(t)) u(t)-r_{2}^{\left(\alpha_{2}\right)}(t)-A_{2}^{\left(\alpha_{2}\right)}(t)+z_{1}(t)=-k_{2} z_{2}(t)
\]
\item Obtain for $u$:
\[
u(t)=-\frac{\mathbf{f}\left(x_{1}(t), x_{2}(t)\right)-r_{2}^{\left(\alpha_{2}\right)}(t)+k_{1}(t) z_{1}^{\left(\alpha_{2}\right)}(t)+z_{1}(t)+k_{2} z_{2}(t)}{\mathbf{g}\left(x_{1}(t), x_{2}(t)\right)}.
\]
\end{enumerate}

In fact, we can generalize the aforementioned procedure for systems of higher order of the form
\[
\begin{aligned}
&x_{i}^{\left(\alpha_{i}\right)}(t)=x_{i+1}(t), \quad\text{for } i=1, \ldots, q-1 \\
&x_{q}^{\left(\alpha_{q}\right)}(t)=\mathbf{f}\left(x_{1}(t), \ldots, x_{q}(t)\right)+\mathbf{g}\left(x_{1}(t), \ldots, x_{q}(t)\right) u(t),
\end{aligned}
\]
where the resulting control law, in this case, is
\begin{equation}\label{eq:control_law_gen}
u(t)=-\frac{\mathbf{f}(x_1(t),\ldots,x_q(t))-r_{q}^{\left(\alpha_{q}\right)}(t)-A_{q}^{\left(\alpha_{q}\right)}(t)+z_{q-1}(t)+k_{q} z_{q}(t)}{\mathbf{g}(x_1(t),\ldots,x_q(t))},
\end{equation}
where $k_{j}>0$, for $j=1, \ldots, q$, and $A_i$ is given by the following recurrence relation: 
\[
\begin{aligned}
A_{1}(t) &=0, \qquad z_{0}(t)=0 \\
A_{i+1}(t) &=-k_{i} z_{i}(t)+A_{i}^{\left(\alpha_{i}\right)}(t)-z_{i-1}(t), \quad\text{for } i=1, \ldots, q-1.
\end{aligned}
\]
Now, the result of applying the control law in detailed in~\eqref{eq:control_law_gen} is 
\begin{equation}\label{eq:gen_cont_law}
\sum_{i=1}^{q} z_{i}(t) z_{i}^{\left(\alpha_{i}\right)}(t)=-\sum_{i=1}^{q} k_{i} z_{i}^{2}(t).
\end{equation}
Finally, to ensure the negativeness of the right-hand side of~\eqref{eq:gen_cont_law} is the same as to ensure the negativity of $\sum_{i=1}^{q} z_{i}(t) \dot{z}_{i}(t)$, and the trajectories in the coordinate system spanned by $z_{1}(t), \ldots, z_{q}(t)$ will converge the origin point.

\section{Adaptive control}\label{sec:adaptive}

\emph{Adaptive control} is the control method used by a controller that must adapt to a controlled system with parameters that either vary over time or are initially uncertain. 
Therefore, it is desirable to have a control law that adapts itself to the changing conditions. 
In other words, adaptive control is a good alternative for industrial applications where the process parameters change, and the controller needs to automatically adapt itself to the new operating conditions. 
This aptitude is called adaptiveness. 
Here, the role of fractional calculus is to design noninteger order adaptation laws or select reference models of noninteger order~\citep{monje2010fractional}. 

\begin{figure}[!ht]
    \centering
    \includegraphics[width=0.6\textwidth]{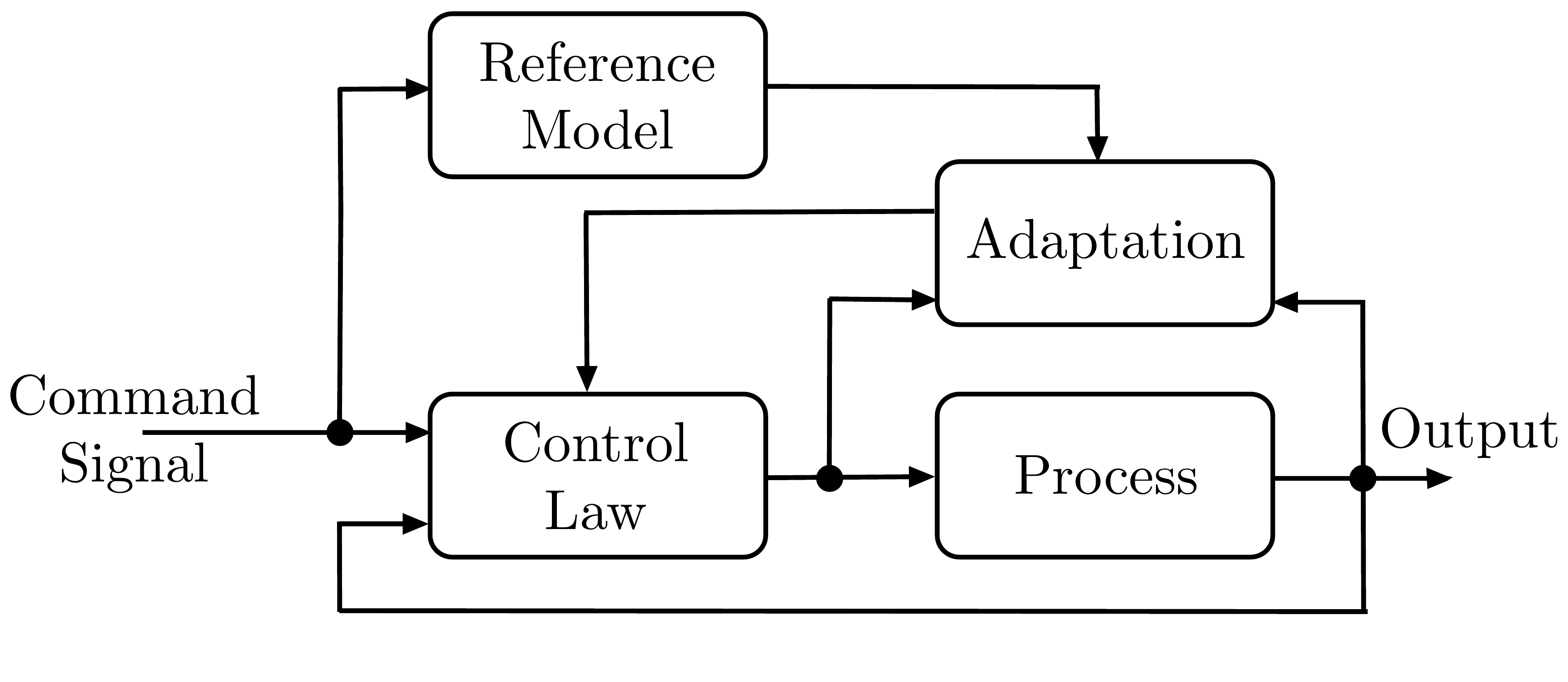}
     \caption{Diagram of the MRAC control scheme.}
    \label{fig:MRAC}
\end{figure}

A broadly adopted adaptive control structure is the so-called \emph{model reference adaptive control} (MRAC) -- see Fig.~\ref{fig:MRAC}. 
This control strategy is grounded on the assumption that the changes in the process parameters are slower than other changes in the closed-loop system. 
Hence, the parameter adjustment mechanism employs the difference between the model output ($y_{m}(t)$) and the process response ($y(t)$) and uses the gradient rule to adjust the parameters of the control law:
\begin{equation}\label{eq:gradient_rule}
\Delta^{\alpha} \phi(t)=-\eta \frac{\partial J(t)}{\partial \phi(t)}=-\eta e(t) \frac{\partial e(t)}{\partial \phi(t)},
\end{equation}
where $\phi(t)$ is a generic parameter of the control law, $e(t)=y(t)-y_{m}(t)$ is the instantaneous model following error, and $J(t)=\frac{e^{2}(t)}{2}$ is the instantaneous performance measure.

It is worth noticing that, when $\alpha=1$ in~\eqref{eq:gradient_rule}, we get the traditional update laws. Additionally, we refer to~\citep{monje2010fractional} for a detailed example considering a  fractional order reference model, where stability is also sought. 
Furthermore, the benefit of utilizing a fractional order setting in MRAC is to achieve a shorter transient regime compared to the classical case. This property might be critical in applications demanding a high-speed response.

\section{System identification}\label{sec:system_identification}

The problem of learning the fractional-order dynamical systems' parameters, i.e., the fractional coefficients and the spatiotemporal matrix, is challenging. Specifically, the determination of the maximum-likelihood poses limitations due to the nonlinearity of the objective. Notwithstanding, some approaches were sucessfully developed in~\citep{gupta2018dealing,gupta2018re,gupta2019learning}, where an approximate solution is based on a variant of the  expectation-maximization algorithm. Nonetheless, such approaches do not enable a finite-time assessment of the uncertainty associated with the parameters, that play a key role in the context of CNS.

Therefore, in what follows, we a recent approach relies on a bilevel iterative bisection scheme~\citep{chatterjee2021learning} to perform identification of the spatial and temporal parameters of a linear discrete-time fractional-order system. First, consider 
\begin{equation}
\label{eq:aug_st_vec}
    \tilde{x}[k] = \begin{bmatrix} x[k] \\ x[k-1] \\ \vdots \\ x[k-p+1] \end{bmatrix}
\end{equation}
as the \emph{augmented} state vector and assuming that the system is \emph{causal}, i.e., the state and disturbances are all considered to be zero before the initial time (i.e., $x[k] = 0$ and $w[k] = 0$ for all $k < 0$), we have
\begin{align}
\label{eq:p_aug_LTI}
    \tilde{x}[k+1] &= \underbrace{\begin{bmatrix} A_0 & \ldots & A_{p-2} & A_{p-1}\\ I & \ldots & 0 & 0\\ \vdots & \ddots & \vdots & \vdots\\ 0 & \ldots & I & 0\end{bmatrix}}_{\tilde{A}}\tilde{x}[k] + \underbrace{\begin{bmatrix}I \\ 0\\ \vdots \\ 0\end{bmatrix}}_{\tilde{B}^w}w[k] \nonumber \\ &= \tilde{A}\tilde{x}[k] + \tilde{B}^w w[k],
\end{align}
for all $k \geq 0$. Note that~\eqref{eq:p_aug_LTI} is an LTI system model, which we refer to as the \emph{$p$-augmented LTI approximation} of~\eqref{eq:fosmain}. 

Having established the $p$-augmented LTI approximation of a DT-FODS in~\eqref{eq:p_aug_LTI}, we can consider the two-level iterative \mbox{bisection-like} approach to identify the spatial and temporal parameters of the \mbox{DT-FODS} in~\eqref{eq:fosmain}. In particular, we start by noting the fact that for the \mbox{Gr\"unwald-Letnikov} definition of the fractional derivative provided in~\eqref{eq:GL_deriv}, $\alpha_i=1$ and $\alpha_i=-1$ can be interpreted, respectively, to be the discretized version of the derivative and the integral for $1 \leq i \leq n$, as defined in the sense of ordinary calculus.

To proceed with a bisection-like approach to identify $\{ \alpha_i \}_{i=1}^n$ and $\tilde{A}$, we first fix the endpoints of the search space for $\alpha_i$ to be $\underline{\alpha_i} = -1$ and $\overline{\alpha_i}=1$ for $1 \leq i \leq n$. We also calculate the value of $\alpha_{c,i} = (\underline{\alpha_i}+\overline{\alpha_i})/2$. Now, given the values of $\underline{\alpha_i}, \overline{\alpha_i}$, and $\alpha_{c,i}$, we calculate, using the ordinary least squares (OLS) technique described in detail below, the row vectors $\underline{\tilde{a}_i}, \overline{\tilde{a_i}}$, and $\tilde{a}_{c,i}$, respectively, that guide the evolution of the states in the $p$-augmented LTI approximation
\begin{equation}
    \tilde{x}_i[k+1] = \tilde{a}_i \tilde{x}_i[k] + \tilde{b}^w_i w_i[k],
\end{equation}
where $\tilde{a}_i = \underline{\tilde{a}_i}$ when $\alpha_i = \underline{\alpha_i}$, $\tilde{a}_i = \overline{\tilde{a}_i}$ when $\alpha_i = \overline{\alpha_i}$, and $\tilde{a}_i = \tilde{a}_{c,i}$ when $\alpha_i = \alpha_{c,i}$ with $\tilde{b}_i^w$ being obtained by extracting the $i$-th row of $\tilde{B}^w$ for $1 \leq i \leq n$.

Next, we propagate the dynamics according to the obtained values of the parameters $\tilde{a}_i$ and calculate the mean squared error (MSE) between the states obtained as a result of the estimated $\tilde{a}_i$'s and the observed states. If the MSE is smaller corresponding to the $\underline{\alpha_i}$ case, then we set $\overline{\alpha_i} = \alpha_{c,i}$. If the MSE is smaller corresponding to the $\overline{\alpha_i}$ case, then we set $\underline{\alpha_i} = \alpha_{c,i}$. This approach is repeated until $\lvert \overline{\alpha_i} - \underline{\alpha_i} \rvert$ does not exceed a certain pre-specified tolerance $\varepsilon$. Algorithm~\ref{alg:sysid} summarizes the procedure of determining the spatial and temporal components of a DT-FODS using the two-level iterative bisection-like approach that we have outlined above.

Therefore, for the estimation of the temporal components of a DT-FODS, we specify the iteration complexity of the bisection-like process and then, we investigate the finite-sample complexity of computing the spatial parameters using a least squares approach.

First, numerical and experimental evidence suggests that the computation of the temporal parameters of a DT-FODS, using, e.g., a wavelet-like technique described in~\citep{flandrin1992wavelet}, does not directly depend on the number of samples or observations used for the aforementioned estimation procedure. Empirical evidence suggests that a small number of samples (usually $30$ to $100$) suffice in order to compute $\{ \alpha_i \}_{i=1}^n$. Furthermore,  we can certify the iteration complexity of the bisection method to find the spatial and temporal parameters of a DT-FODS. Specifically, the bisection-based technique detailed above to find the temporal components of a DT-FODS is minmax optimal and the number $\nu$ of iterations needed in order to achieve a certain specified tolerance $\varepsilon$ when this technique is used is bounded above by
\begin{equation}
    \nu \leq \left\lceil \log_2 \left( \frac{2}{\varepsilon} \right) \right\rceil.
\end{equation}

\begin{algorithm}
\caption{Learning the parameters of a DT-FODS}
\begin{algorithmic}[1]
\label{alg:sysid}
\FOR{$i=1$ to $n$}
\STATE{Initialize $\underline{\alpha_i}=-1$, $\overline{\alpha_i}=1$, and tolerance $\varepsilon$.}
\STATE{Calculate $\alpha_{c,i} = (\underline{\alpha_i}+\overline{\alpha_i})/2$.}
\STATE{Given the above values of $\underline{\alpha_i}$, $\overline{\alpha_i}$, and $\alpha_{c,i}$, find, using the ordinary least squares (OLS) method, the row vectors $\underline{\tilde{a}_i}, \overline{\tilde{a_i}}$, and $\tilde{a}_{c,i}$, respectively, that guide the evolution of the states in the $p$-augmented LTI approximation $\tilde{x}_i[k+1] = \tilde{a}_i \tilde{x}_i[k] + \tilde{b}^w_i w_i[k]$.}
\STATE{Propagate the dynamics according to the obtained OLS estimates and calculate the mean squared error (MSE) between the propagated states and the observed state trajectory.}
\IF{MSE is lesser for the $\underline{\alpha_i}$ case}
\STATE{Set $\overline{\alpha_i} = \alpha_{c,i}$.}
\ELSIF{MSE is lesser for the $\overline{\alpha_i}$ case}
\STATE{Set $\underline{\alpha_i} = \alpha_{c,i}$.}
\ENDIF
\STATE{Terminate if $\lvert \overline{\alpha_i} - \underline{\alpha_i} \rvert < \varepsilon$, else return to step 3.}
\ENDFOR
\end{algorithmic}
\end{algorithm}


Secondly, we can now delve into the problem of identifying the spatial parameters using a least squares-like approach and its finite-time guarantees. We start with the $p$-augmented LTI model of~\eqref{eq:p_aug_LTI}, i.e.,
\begin{equation}
    \tilde{x}[k+1] = \tilde{A}\tilde{x}[k] + \tilde{B}^w w[k].
\end{equation}
The OLS method then outputs the matrix $\underline{\tilde{A}}[K]$ as the solution of the following optimization problem
\begin{equation}
    \underline{\tilde{A}}[K] \coloneqq \argmin_{\tilde{A} \in \mathbb{R}^{d \times d}} \sum_{k=1}^K \frac{1}{2} \left\| \tilde{x}[k+1] - \tilde{A} \tilde{x}[k] \right\|_2^2,
\end{equation}
by observing the state trajectory of~\eqref{eq:p_aug_LTI}, i.e., $\{ x[0], x[1], \ldots, x[K+1] \}$, and the process noise $w[k]$ being independent and identically distributed (i.i.d.) zero-mean Gaussian.


Thus, prior to characterizing the sample complexity of the OLS method for the $p$-augmented LTI approximation of the \mbox{DT-FODS}, we define a few quantities of interest.
The \emph{finite-time controllability Gramian} of the approximated system~\eqref{eq:p_aug_LTI}, $W_t$, is defined by
\begin{equation}
    W_t \coloneqq \sum_{j=0}^{t-1} \tilde{A}^j (\tilde{A}^j)^{\mathsf{T}}.
\end{equation}
Intuitively, the controllability Gramian gives a quantitative measure of how much the system is excited when induced by the process noise $w[k]$ acting as an input to the system.

Additionally, given a symmetric matrix $A \in \mathbb{R}^{d \times d}$, we define $\lambda_{\max}(A)$ and $\lambda_{\min}(A)$ to denote, respectively, the maximum and minimum eigenvalues of the matrix $A$.

Lastly, for any square matrix $A \in \mathbb{R}^{d \times d}$, the \emph{spectral radius} of the matrix $A$, $\rho(A)$, is given by the largest absolute value of its eigenvalues. Also, the \emph{operator norm} of a matrix is denoted by 
$${\displaystyle \|A\|_{\text{op}}=\inf\{c\geq 0:\|Av\|\leq c\|v\|{\mbox{ for all }}v\in V\}.}$$



Hence, we have the following result that characterizes the sample complexity of the above OLS method for the DT-FODS approximation.
Fix $\delta \in (0,1/2)$ and consider the \mbox{$p$-augmented} system in~\eqref{eq:p_aug_LTI}, where $\tilde{A} \in \mathbb{R}^{d \times d}$ is a marginally stable matrix (i.e., $\rho(\tilde{A}) \leq 1$) and $w[k] \thicksim \mathcal{N}(0,\sigma^2 I)$. Then, there exist universal constants $c,C > 0$ such that,
\begin{align}
\label{eq:bound_without_inputs}
\mathbb{P}\Bigg[ \left\| \underline{\tilde{A}}[K] - \tilde{A} \right\|_{\mathrm{op}} \leq  \frac{C}{\sqrt{K \lambda_{\min }\left(W_{k}\right)}} \nonumber  \times &\sqrt{d \log \left( \frac{d}{\delta} \right) +\log \operatorname{det}\left(W_{K} W_{k}^{-1}\right)} \Bigg]& \\  & \geq 1 - \delta,
\end{align}
for any $k$, such that
\begin{equation}
    \frac{K}{k} \geq c\left(d \log \left( \frac{d}{\delta} \right) + \log \operatorname{det}\left(W_{K} W_{k}^{-1}\right)\right)
\end{equation}
holds.

\begin{remark}
We note here that although the operator norm parameter estimation error in~\eqref{eq:bound_without_inputs} is stated in terms of $\tilde{A}$, the operator norm errors, associated with the matrices $A_0, A_1, \ldots, A_{p-1}$, are strictly lesser compared to $\left\| \underline{\tilde{A}}[K] - \tilde{A} \right\|_{\mathrm{op}}$, since $A_0, A_1, \ldots, A_{p-1}$ are submatrices of $\tilde{A}$, and for any operator norm, the operator norm of a submatrix is upper bounded by one of the whole matrix (see Lemma A.9 of~\citep{foucart2013} for a proof).
\end{remark}

Additionally, it is worth mentioning that a similar finite-sample complexity bound similar to the one presented before can also be derived when we consider the ordinary least squares identification of the spatial parameters of a DT-FODS with inputs. For instance, within the purview of epileptic seizure mitigation using intracranial EEG data, the objective of the former is to suppress the overall length or duration of an epileptic seizure. Thus, the goal is steering the state of the neurophysiological system in consideration away from seizure-like activity, using a control strategy like model predictive control~\citep{chatterjee2020fractional}.

\section{Minimum-energy state estimation}\label{sec:state_estimation}

Most of the estimators that exist for fractional-order dynamical systems are obtained under the assumption that the disturbance and noise has Gaussian distribution~\citep{sabatier2012observability,sierociuk2006fractional,Safari1,Safari2,miljkovic2017ecg,najar2009discrete}. Meanwhile, such assumption is not realistic in the context of neural systems as disturbance frequencies can only lie within a specific frequency band. Therefore, in what follows, we present the so-called \emph{minimum-energy state estimation}, where it is assumed that the disturbance and noise are unknown, but deterministic and bounded uncertainties.

Now, consider a left-bounded sequence $\{ x[k] \}_{k \in \mathbb{Z}}$ over $k$, i.e., with $\limsup\limits_{k \to -\infty} \| x[k] \| < \infty$. Then, the \mbox{Gr\"{u}nwald-Letnikov} fractional-order difference, for any $\alpha \in \mathbb{R}^{+}$, can be re-written  as
\begin{equation}
\label{eq:diff_op}
\begin{split}
\Delta^{\alpha} x[k] \coloneqq \sum_{j=0}^{\infty} c_j^{\alpha} x[k-j],\quad c_j^{\alpha} = (-1)^j \binom{\alpha}{j}, \quad \\ \binom{\alpha}{j} = \begin{cases} 1 &\mbox{if } j = 0, \\
                    \displaystyle\prod_{i=0}^{j-1} \frac{\alpha-i}{i+1} = \frac{\Gamma (\alpha+1)}{\Gamma (j+1) \Gamma (\alpha-j+1)} & \mbox{if } j > 0, \end{cases}
\end{split}
\end{equation}
for all $j \in \mathbb{N}$. The summation in~\eqref{eq:diff_op} is well-defined from the uniform boundedness of the sequence $\{ x[k] \}_{k \in \mathbb{Z}}$ and the fact that $|c^{\alpha}_j| \leq \frac{\alpha^j}{j!}$, which implies that the sequence $\{ c^{\alpha}_j \}_{j \in \mathbb{N}}$ is absolutely summable for any $\alpha \in \mathbb{R}^{+}$~\citep{alessandretti2020finite,sopasakis2017stabilising}.

With the above ingredients, a {discrete-time fractional-order dynamical network with additive disturbance} can be described, respectively, by the state evolution and output equations
\begin{subequations}
\begin{equation}
\label{eq:fos_model}
\sum_{i=1}^{l} A_{i} \Delta^{a_{i}} x[k+1]=\sum_{i=1}^{r} B_{i} \Delta^{b_{i}} u[k]+\sum_{i=1}^{s} G_{i} \Delta^{g_{i}} w[k],
\end{equation}
\begin{equation}
\label{eq:fos_model_op}
z[k] = C'_k x[k] + v'[k],
\end{equation}
\end{subequations}
with the variables $x[k] \in \mathbb{R}^n$, $u[k] \in \mathbb{R}^m$, and $w[k] \in \mathbb{R}^p$ denoting the {state}, {input}, and {disturbance} vectors at time step $k \in \mathbb{N}$, respectively. The scalars $a_i \in \mathbb{R}^{+}$ with $1 \leq i \leq l$, $b_i \in \mathbb{R}^{+}$ with $1 \leq i \leq r$, and $g_i \in \mathbb{R}^{+}$ with $1 \leq i \leq s$ are the {\mbox{fractional-order} coefficients} corresponding, respectively, to the state, the input, and the disturbance. The vectors $z[k], v'[k] \in \mathbb{R}^q$ denote, respectively, the {output} and {measurement disturbance} at time step $k \in \mathbb{N}$. We assume that the (unknown but deterministic) disturbance vectors are bounded as
\begin{equation}
    \| w[k] \| \leq b_w, \| v'[k] \| \leq b_{v'}, \; k \in \mathbb{N},
\end{equation}
for some scalars $b_w, b_{v'} \in \mathbb{R}^{+}$. We also assume that the control input $u[k]$ is known for all time steps $k \in \mathbb{N}$. We denote by $x[0] = x(0)$ the initial condition of the state at time $k = 0$. In the computation of the fractional-order difference, we assume that the system is \emph{causal}, i.e., the state, input, and disturbances are all considered to be zero before the initial time (i.e., $x[k] = 0, u[k] = 0$, and $w[k] = 0$ for all $k < 0$).

Next, consider the quadratic weighted least-squares objective function
\begin{align}
\displaystyle\mathcal{J}\left( x[0], \{ w[i] \}_{i=0}^{N-1}, \{ v'[j] \}_{j=1}^N \right) &= \sum_{i=0}^{N-1} w[i]^{\mathsf{T}} Q_i^{-1} w[i] + \sum_{j=1}^{N} v'[j]^{\mathsf{T}} R_j^{-1} v'[j] \nonumber \\ &+ (x[0] - \hat{x}_0)^{\mathsf{T}} P_0^{-1} (x[0] - \hat{x}_0),
\end{align}
subject to the constraints
\begin{subequations}
\begin{equation}
\label{eq:fos_model_main1}
\sum_{i=1}^{l} A_{i} \Delta^{a_{i}} x[k+1]=\sum_{i=1}^{r} B_{i} \Delta^{b_{i}} u[k]+\sum_{i=1}^{s} G_{i} \Delta^{g_{i}} w[k]
\end{equation}
and
\begin{equation}
\label{eq:fos_model_main2}
z[k] = C'_k x[k] + v'[k],
\end{equation}
\end{subequations}
for some $N \in \mathbb{N}$, with the weighting matrices $Q_i$ $(0 \leq i \leq N-1), R_j$ $(1 \leq j \leq N)$, and $P_0$ chosen to be symmetric and positive definite, and $\hat{x}_0$ chosen to be the \emph{a priori} estimate of the system's initial state. The minimum-energy estimation procedure seeks to solve the following optimization problem
\begin{mini}|l|
  {\scriptstyle \{ x[k] \}_{k=0}^{N}, \{ w[i] \}_{i=0}^{N-1}, \{ v'[j] \}_{j=1}^{N}}{\mathcal{J}\left( x[0], \{ w[i] \}_{i=0}^{N-1}, \{ v'[j] \}_{j=1}^N \right)}{}{}
  \addConstraint{\eqref{eq:fos_model_main1} \: \mathrm{and} \: \eqref{eq:fos_model_main2}},
\label{eq:opt_prob_main1}
\end{mini}
for some $N \in \mathbb{N}$.

To derive the solution to~\eqref{eq:opt_prob_main1}, we first start with some alternative formulations of the discrete-time fractional-order dynamical network (DT-FODN) in~\eqref{eq:fos_model} and~\eqref{eq:fos_model_op} and relevant definitions that will be used in the sequel. Then, we present the solution  and  some additional properties of the derived solution, i.e., the exponential input-to-state stability of the estimation error.

In what follows, we consider the  mild technical assumption that  $\displaystyle\sum_{i=1}^l A_i$ is invertible.
Additionally, we consider a truncation of the last $\mathfrak{v}$ temporal components of~\eqref{eq:fos_model}, which we will refer to as the $\mathfrak{v}$-approximation for the DT-FODN. That being said, the DT-FODN model in~\eqref{eq:fos_model} can be equivalently written as
\begin{equation}
    x[k+1]=\sum_{j=1}^{\infty} \check{A}_j x[k-j+1] + \sum_{j=0}^{\infty} \check{B}_j u[k-j] + \sum_{j=0}^{\infty} \check{G}_j w[k-j],
\end{equation}
where $\check{A}_j = -\hat{A}_0^{-1} \hat{A}_j$, $\check{B}_j = \hat{A}_0^{-1} \hat{B}_j$, and $\check{G}_j = \hat{A}_0^{-1} \hat{G}_j$ with $\hat{A}_j = \sum_{i=1}^l A_i c_j^{a_i}$, $\hat{B}_j = \sum_{i=1}^r B_i c_j^{b_i}$, and \mbox{$\hat{G}_j = \sum_{i=1}^s G_i c_j^{g_i}$}. Furthermore, for any positive integer $\mathfrak{v} \in \mathbb{N}^{+}$, the DT-FODN model in~\eqref{eq:fos_model} can be recast as
\begin{subequations}
\label{eq:sys}
\begin{equation}
\label{eq:v_app1}
    \tilde{x}[k+1] = \tilde{A}_{\mathfrak{v}} \tilde{x}[k] + \tilde{B}_{\mathfrak{v}} u[k] + \tilde{G}_{\mathfrak{v}} r[k], \qquad \tilde{x}[0] = \tilde{x}_0,
\end{equation}
\begin{equation}
\label{eq:v_app3}
    y[k+1] = C_{k+1} \tilde{x}[k+1] + v[k+1],
\end{equation}
\end{subequations}
where
\begin{equation}
\label{eq:v_app2}
    r[k] = \sum_{j={\mathfrak{v}}+1}^{\infty} \check{A}_j x[k-j+1] + \sum_{j={\mathfrak{v}}+1}^{\infty} \check{B}_j u[k-j] + \sum_{j=0}^{\infty} \check{G}_j w[k-j],
\end{equation}
with the augmented state vector $\tilde{x}[k] = [ x[k]^{\mathsf{T}}, \ldots, x[k-\mathfrak{v}+1]^{\mathsf{T}}, u[k-1]^{\mathsf{T}},\ldots,u[k-\mathfrak{v}]^{\mathsf{T}} ]^{\mathsf{T}} \in \mathbb{R}^{\mathfrak{v} \times (n+m)}$ and appropriate matrices $\tilde{A}_{\mathfrak{v}}, \tilde{B}_{\mathfrak{v}}$, and $\tilde{G}_{\mathfrak{v}}$, where $\tilde{x}_0 = [x_0^{\mathsf{T}},0,\ldots,0]^{\mathsf{T}}$ denotes the initial condition. The matrices $\tilde{A}_{\mathfrak{v}}$ and $\tilde{B}_{\mathfrak{v}}$ are formed using the terms $\{ \check{A}_j \}_{1 \leq j \leq \mathfrak{v}}$ and $\{ \check{B}_j \}_{1 \leq j \leq \mathfrak{v}}$, while the remaining terms $\{ \check{G}_j \}_{1 \leq j < \infty}$ and the state and input components not included in $\tilde{x}[k]$ are absorbed into the term $\tilde{G}_{\mathfrak{v}} r[k]$. Furthermore, we refer to~\eqref{eq:v_app1} as the \emph{$\mathfrak{v}$-approximation} of the DT-FODN presented in~\eqref{eq:fos_model}.


To obtain the minimum-energy estimator, let us consider the quadratic weighted least-squares objective function
\begin{align}
\label{eq:objective}
\mathcal{J}\left( \tilde{x}[0], \{ r[i] \}_{i=0}^{N-1}, \{ v[j] \}_{j=1}^N \right) &= \sum_{i=0}^{N-1} r[i]^{\mathsf{T}} Q_i^{-1} r[i] + \sum_{j=1}^{N} v[j]^{\mathsf{T}} R_j^{-1} v[j] \nonumber \\ &+ (\tilde{x}[0] - \hat{x}_0)^{\mathsf{T}} P_0^{-1} (\tilde{x}[0] - \hat{x}_0),
\end{align}
subject to the constraints
\begin{subequations}
\label{eq:syscon}
\begin{align}
\begin{split}
    \label{eq:syscon1}
    \bar{x}[k+1] &= \tilde{A}_{\mathfrak{v}} \bar{x}[k] + \tilde{B}_{\mathfrak{v}} u[k] + \tilde{G}_{\mathfrak{v}} \bar{r}[k],    
\end{split}\\
\begin{split}
    \label{eq:syscon2}
    y[k+1] &= C_{k+1} \bar{x}[k+1] + \bar{v}[k+1],
\end{split}
\end{align}
\end{subequations}
for some $N \in \mathbb{N}$. The weighting matrices $Q_i$ $(0 \leq i \leq N-1)$ and $R_j$ $(1 \leq j \leq N)$ are chosen to be symmetric and positive definite. The term $\hat{x}_0$ denotes the \emph{a priori} estimate of the (unknown) initial state of the system, with the matrix $P_0$ being symmetric and positive definite.

Subsequently, we  consider the weighted \mbox{least-squares} optimization problem
\begin{mini}|l|
  {\scriptstyle \{ \bar{x}[k] \}_{k=0}^{N}, \{ \bar{r}[i] \}_{i=0}^{N-1}, \{ \bar{v}[j] \}_{j=1}^{N}}{\mathcal{J}\left( \tilde{x}[0], \{ r[i] \}_{i=0}^{N-1}, \{ v[j] \}_{j=1}^N \right)}{}{}
  \addConstraint{\eqref{eq:syscon1} \: \text{and} \: \eqref{eq:syscon2}},
\label{eq:opt_prob}
\end{mini}
for some $N \in \mathbb{N}$. 
Denote by $\hat{x}[k]$ the state vector that corresponds to the solution of the optimization problem~\eqref{eq:opt_prob}. Then, $\hat{x}[k]$ satisfies the recursion
\begin{equation}
\label{eq:opt_sol}
    \hat{x}[k+1] = \tilde{A}_{\mathfrak{v}} \hat{x}[k] + \tilde{B}_{\mathfrak{v}} u[k] + K_{k+1} \left( y[k+1] - C_{k+1} \left( \tilde{A}_{\mathfrak{v}} \hat{x}[k] + \tilde{B}_{\mathfrak{v}} u[k] \right) \right),
\end{equation}
given $0 \leq k \leq N-1$, with initial conditions specified for $\hat{x}_0$ and $\{ u[j] \}_{j=0}^k$, and with the update equations
\begin{subequations}
\label{eq:filter_updates}
\begin{equation}
    K_{k+1} = M_{k+1} C_{k+1}^{\mathsf{T}} ( C_{k+1} M_{k+1} C_{k+1}^{\mathsf{T}} + R_{k+1} )^{-1},
\end{equation}
\begin{equation}
    M_{k+1} = \tilde{A}_{\mathfrak{v}} P_k \tilde{A}_{\mathfrak{v}}^{\mathsf{T}} + \tilde{G}_{\mathfrak{v}} Q_k \tilde{G}_{\mathfrak{v}}^{\mathsf{T}},
\end{equation}
and
\begin{align}
\label{eq:P_update}
    P_{k+1} &= (I - K_{k+1} C_{k+1}) M_{k+1} (I - K_{k+1} C_{k+1})^{\mathsf{T}} +  K_{k+1} R_{k+1} K_{k+1}^{\mathsf{T}} \nonumber \\ &= (I - K_{k+1} C_{k+1}) M_{k+1},
\end{align}
\end{subequations}
with symmetric and positive definite $P_0$.

Notice that the dynamics of the recursion in~\eqref{eq:opt_sol} (with the initial conditions on $\hat{x}_0$ and the values of $\{ u[j] \}_{j=0}^k$ being known) along with the update equations~\eqref{eq:filter_updates} together solve~\eqref{eq:opt_prob}. It is interesting to note here that the output term $y[k+1]$ presented in~\eqref{eq:syscon2} and~\eqref{eq:opt_sol} is the output of the $\mathfrak{v}$-approximated system~\eqref{eq:sys}, which, in turn, is simply a subset of the outputs $z[k+1]$ obtained from~\eqref{eq:fos_model_op}, truncated $\mathfrak{v}$ time steps in the past, provided $v[k]$ and $C_k$ are formed from the appropriate blocks of $v'[k]$ and $C'_k$ for all $k \in \mathbb{N}$.

Secondly, the minimum-energy estimator has exponential input-to-state stability of the estimation error.



In order to prove the exponential input-to-state stability of the minimum-energy estimation error, we need to consider the following mild technical assumptions. Specifically, there exist constants $\underline{\alpha}, \overline{\alpha}, \beta, \gamma \in \mathbb{R}^{+}$ such that
\begin{equation}
    \underline{\alpha} I \preceq \tilde{A}_{\mathfrak{v}} \tilde{A}_{\mathfrak{v}}^{\mathsf{T}} \preceq \overline{\alpha} I, \quad \tilde{G}_{\mathfrak{v}} \tilde{G}_{\mathfrak{v}}^{\mathsf{T}} \preceq \beta I, \:\: \text{and} \quad C_k^{\mathsf{T}} C_k \preceq \gamma I,
\end{equation}
for all $k \in \mathbb{N}$.

Additionally, notice that the \emph{state transition matrix} for the dynamics in~\eqref{eq:v_app1} is given by
\begin{equation}
    \Phi(k,k_0) = \tilde{A}_{\mathfrak{v}}^{(k-k_0)}, \quad \text{with} \quad \Phi(k_0,k_0) = I,
\end{equation}
for all $k \geq k_0 \geq 0$. We also consider the \emph{discrete-time controllability Gramian} associated with the dynamics~\eqref{eq:v_app1} described by
\begin{equation}
\label{eq:contr_gramian}
    W_c(k,k_0) = \sum_{i = k_0}^{k-1} \Phi (k,i+1) \tilde{G}_{\mathfrak{v}} \tilde{G}_{\mathfrak{v}}^{\mathsf{T}} \Phi^{\mathsf{T}} (k,i+1),
\end{equation}
and the \emph{discrete-time observability Gramian} associated with~\eqref{eq:v_app1} to be
\begin{equation}
\label{eq:obsv_gramian}
    W_o(k,k_0) = \sum_{i = k_0+1}^{k} \Phi^{\mathsf{T}} (i,k_0) C_i^{\mathsf{T}} C_i \Phi (i,k_0),
\end{equation}
for $k \geq k_0 \geq 0$. We also make the following assumptions regarding \emph{complete uniform controllability} and \emph{complete uniform observability} of the \mbox{$\mathfrak{v}$-approximated} system in~\eqref{eq:v_app1}.

As such, we have to also consider that the $\mathfrak{v}$-approximated system~\eqref{eq:v_app1} is completely uniformly controllable, i.e., there exist constants $\delta \in \mathbb{R}^{+}$ and $N_c \in \mathbb{N}^{+}$ such that
\begin{equation}
    W_c (k+N_c,k) \succeq \delta I,
\end{equation}
for all $k \geq 0$. And, similarly, the $\mathfrak{v}$-approximated system~\eqref{eq:v_app1} is completely uniformly observable, i.e., there exist constants \mbox{$\varepsilon \in \mathbb{R}^{+}$} and $N_o \in \mathbb{N}^{+}$ such that
\begin{equation}
    W_o (k+N_o,k) \succeq \varepsilon \Phi^{\mathsf{T}} (k+N_o,k) \Phi(k+N_o,k),
\end{equation}
for all $k \geq 0$.

Next, we also present an assumption certifying lower and upper bounds on the weight matrices $Q_k$ and $R_{k+1}$ in~\eqref{eq:objective}. That is,  without loss of generality, we assume that the weight matrices $Q_k$ and $R_{k+1}$ satisfy
\begin{equation}
    \underline{\vartheta} I \preceq Q_k \preceq \overline{\vartheta} I \quad \text{and} \quad \underline{\rho} I \preceq R_{k+1} \preceq \overline{\rho} I,
\end{equation}
for all $k \geq 0$ and constants $\underline{\vartheta}, \overline{\vartheta}, \underline{\rho}, \overline{\rho} \in \mathbb{R}^+$.

Hence, it is possible to we establish lower and upper bounds for the matrix $P_k$ required to show that the estimation error is exponentially input-to-state stable. Specifically,  the minimum-energy estimation error $e[k]$, given by
\begin{equation}
\label{eq:err_defn}
    e[k] = \hat{x}[k] - \tilde{x}[k],
\end{equation}
is such that there exist constants $\sigma, \tau, \chi, \psi \in \mathbb{R}^{+}$ with $\tau < 1$ such that the estimation error $e[k]$ satisfies
\begin{equation}
\label{eq:ISS_theorem}
    \| e[k] \| \leq \max \Bigg\{ \sigma \tau ^{k-k_0} \| e[k_0] \|,\: \chi \max_{k_o \leq i \leq k-1} \| r[i] \|,\: \psi \max_{k_o \leq j \leq k-1} \| v[j+1] \| \Bigg\}
\end{equation}
for all $k \geq k_0 \geq \max\{ N_c, N_o \}$.

It is interesting to note that the bound on the estimation error $e[k]$ in~\eqref{eq:ISS_theorem} actually depends on $\| r[i] \|$, where $k_0 \leq i \leq k-1$ for all $i \in \mathbb{N}$. In fact, a distinguishing feature of DT-FODN is the presence of a finite non-zero disturbance term in the \mbox{input-to-state} stability bound of the tracking error when tracking a state other than the origin. This disturbance is dependent on the upper bounds on the non-zero reference state being tracked as well as the input. While the linearity of the \mbox{Gr\"unwald-Letnikov} fractional-order difference operator allows one to mitigate this issue in the case of tracking a non-zero exogenous state by a suitable change of state and input coordinates, this approach is not one we can pursue in this paper, since the state we wish to estimate is unknown. However, it can be shown that as the value of $\mathfrak{v}$ in the $\mathfrak{v}$-approximation increases, the upper bound associated with $\| r[i] \|$ decreases drastically since the $\mathfrak{v}$-approximation gives us progressively better representations of the unapproximated system. This further implies that $\| r[i] \|$ in~\eqref{eq:ISS_theorem} stays bounded, with progressively smaller upper bounds associated with $\| r[i] \|$ (and hence, $\| e[k] \|$) with increasing $\mathfrak{v}$.

Lastly, the estimation error associated with the minimum-energy estimation process in~\eqref{eq:err_defn} is defined in terms of the state of the $\mathfrak{v}$-approximated system $\tilde{x}[k]$. In reality, as detailed above, with larger values of $\mathfrak{v}$, the $\mathfrak{v}$-approximated system approaches the real network dynamics, and thus we obtain an expression for the estimation error with respect to the real system in the limiting case, where the input-to-state stability bound as presented in~\eqref{eq:ISS_theorem} holds.

\section{Fractional optimal control}\label{sec:fractional_optimal_control}
Fractional optimal control finds the optimal control strategy to manipulate a fractional-order dynamical system to achieve a specific goal. Usually the goal is to achieve a certain desired state behavior while minimizing the amount of control effort~\citep{riewe1996nonconservative}. The fractional optimal control problem with a finite-time horizon can be formulated as follows: 

\begin{equation}
{\small
\begin{aligned}
& \text{(cost function)} & &\underset{\mathbf{u}\quad\,}{\text{minimize\hspace{0.5cm}}}
& & \int_{t_0}^{T} (\mathbf{x}(t)-\mathbf{x}_{d}(t))^{\intercal}Q(\mathbf{x}(t)-\mathbf{x}_{d}(t))+\mathbf{u}(t)^{\intercal}R\mathbf{u}(t) \mathrm{~d}t\\
& \text{(constraints)} & &\text{subject to}
&  &\Delta^\alpha \mathbf{x}(t) = A\mathbf{x}(t) + B\mathbf{u}(t) \\
& & & & & \text{other linear constraints on } \mathbf{x}(t) \text{ and } \mathbf{u}(t),
\end{aligned}
}
\end{equation}
where $\mathbf{x}(t)\in\mathbb{R}^{n}$ is the state of the system, $\mathbf{x}_{d}(t)\in\mathbb{R}^{n}$ is the desired state of the system, $\mathbf{u}(t)\in\mathbb{R}^{m}$ is the control input, $Q$ is the cost on the state achieving the desired behavior, $R$ is the cost on the control effort, $\Delta^\alpha$ is the Caputo fractional-order derivative, $A\in\mathbb{R}^{n\times n}$ is the state matrix, and $B\in\mathbb{R}^{n\times m}$ is the control input  matrix.

Many mathematical techniques for solving fractional optimal control problems have been proposed, including numerical solvers~\citep{nemati2019numerical,baleanu2009central,agrawal2007hamiltonian,agrawal2004general} and discrete methods~\citep{almeida2015discrete}. Other works have considered fractional optimal control using the following schemes, including distributed fractional optimal control~\citep{zaky2017formulation}, finite-time horizon~\citep{biswas2011fractional}, multi-dimensional~\citep{agrawal2010fractional}, Euler-Lagrange formulation~\citep{frederico2008fractional,frederico2007formulation,torres2012introduction,agrawal2002formulation}, and reinforcement learning~\citep{gupta2021non}. 
Furthermore, fractional optimal control has been used in the followoing applications cloud computing~\citep{ghorbani2014prediction}, regulating diabetes~\citep{ghorbani2014reducing,ghorbani2013cyber}, cyber-physical systems~\citep{bogdan2011towards}, regulating heart disease~\citep{bogdan2012implantable,bogdan2013pacemaker}, data-centers-on-chip~\citep{bogdan2015mathematical}, and power management~\citep{bogdan2012optimal,bogdan2013dynamic}, and chemical processing plants~\citep{petravs2021novel}.

Fractional optimal control is at the core of receding horizon approaches referred to as model predictive control, and overviewed in more detail next.

\section{Model predictive control}\label{sec:model_predictive_control}


Model predictive control (MPC) is a control strategy that allows the control of processes while satisfying a set of constraints.  At its core, MPC uses explicit process models (which may be linear or nonlinear) to predict how a plant will respond to arbitrary inputs. For each instant of time, an MPC algorithm seeks to optimize plant behavior in the future by computing a series of control inputs over a time horizon called the \emph{prediction horizon} by solving an optimization problem -- often with constraints. Once this step is complete, the computed control inputs corresponding to the first subsection of the prediction horizon (called the \emph{control horizon}) are then sent to the plant. This procedure is then repeated at subsequent control intervals~\citep{qin2003survey}. This receding horizon strategy implicitly introduces \emph{closed-loop feedback}.

Next,  we consider the case where the  predictive model  is a linear fractional-order system. Based on the state signal's evolution predicted by the model, and by regarding the impact of an arbitrary control input signal in the state's evolution, we can set out to adapt the stimulation signal in real-time by choosing the parameters that lead to stimulation signals within a safe range towards optimizing some measure of performance that encapsulates the goal of steering abnormal activity to normal ranges. In general, however, our predictive model will not precisely match the real dynamics of the system. Therefore, our proposed stimulation strategy will periodically re-evaluate the current estimated state and corresponding predictions and re-compute the appropriate optimal stimulation strategy. 

First, in the fractional-order model predictive control framework, we will focus on the design of a model predictive controller for a (possibly time-varying) discrete-time fractional-order dynamical system model
\begin{equation}
    \Delta^\alpha x[k+1] = A_k x[k] + B_k u[k] + B^w_k w[k],
    \label{eq:fosmodel}
\end{equation}
where $w[k]$ denotes a sequence of independent and identically distributed random vectors, following an $\mathcal{N}(0,\Sigma)$ distribution (with the covariance matrix $\Sigma \in \mathbb{R}^{n \times n}$) and $B^w_k$ denotes the matrix of weights that scales the noise term $w[k]$. The objective is to design the feedback controller such that it minimizes a quadratic cost functional of the input and state vectors over a finite time horizon $P$ (the prediction horizon). In other words, the objective is to determine the sequence of control inputs $u[k],\ldots,u[k+P-1]$ that minimizes a quadratic cost function of the form
\begin{equation}
{\small
\begin{aligned}
& \text{(cost function)} & &\underset{u[k],\ldots,u[k+P-1]\quad\,}{\text{minimize\hspace{0.5cm}}}
& &\mathbb{E}\Big\{ \sum_{j=1}^{P}\|x[k+j] \|_{{Q}_{k+j}}^2 + \sum_{j=1}^P{c}_{k+j}^{\mathsf{T}}x[k+j] \\ 
& & & & &+ \sum_{j=0}^{P-1}\| u[k+j] \|_{{R}_{k+j}}^2 \Big\}\\
& \text{(constraints)} & &\text{subject to}
&  &x[k] = \text{ observed or estimated current state}\\
& & & & & \Delta^\alpha x[k+j+1] = A_{k+j} x[k+j] + B_{k+j} u[k+j] \\
& & & & & + B^w_{k+j} w[k+j],\\ 
& & & & & j=0,1,\ldots,P-1,\\
& & & & & \text{other linear constraints on } x[k+1],\ldots,x[k+P],\\
& & & & & u[k],\ldots,u[k+P-1],
\end{aligned}
}
\label{eq:mpc123_jne}
\end{equation}
where ${Q}_{k+1},\ldots,{Q}_{k+P}\in\mathbb{R}^{n\times n}$ and ${R}_{k},\ldots,{R}_{k+P-1}\in\mathbb{R}^{n_u\times n_u}$ are given positive semidefinite matrices. Here, \mbox{$Q\in\mathbb{R}^{n\times n}$} is a \emph{positive semidefinite} matrix if $x^{\mathsf{T}} Qx\geq 0$, for every $x\in\mathbb{R}^n$, and $\|x\|_Q = \sqrt{x^{\mathsf{T}} Qx}$ in that case. 


The quadratic term on the input, which represents the electrical neurostimulation signal, is intended to add a penalization term for stimulating the patient too harshly, since this may be unsafe, create discomfort for the patient, or result in harmful psychological effects~\citep{moratti2014adverse}. It is also interesting to note that even if we need the estimation of the system states in the above problem, the presence of a separation principle for \mbox{discrete-time} \mbox{fractional-order} systems~\citep{chatterjee2019sep} gives us guarantees that we can perform model predictive control with state estimation for these systems.

Note that, here, $P$ is called the \emph{prediction horizon}, and the framework only deploys the control strategy associated with the first $M$ time steps (referred to as the \emph{control horizon}). Simply speaking, after we reach state $x[k+M-1]$, we update $k$ with $k+M-1$ and recompute the new solution. This way, we have robust solutions, since, by design, the optimal strategy is constantly being re-evaluated based on the short-term control action implementation of a long-term prediction~\citep{Bequette2013,petravs2021novel}.

\section{Applications in cyber-neural systems}\label{sec:applications}

\subsection{System identification}

\begin{figure}
    \centering
    \includegraphics[width=0.9\textwidth]{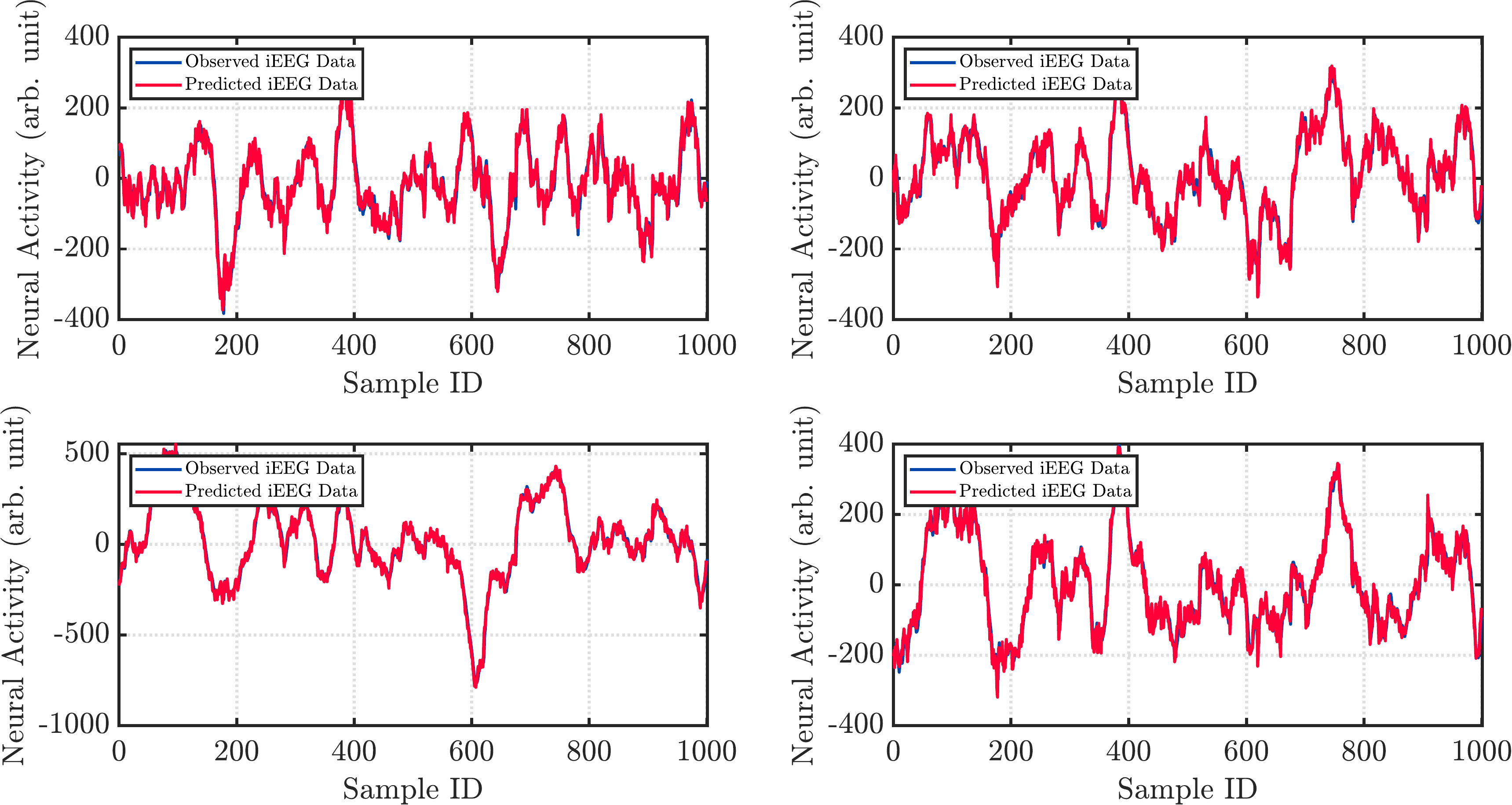}
    \caption{Performance of our system identification approach on real-life intracranial EEG data.}
    \label{fig:eeg_sysid}
\end{figure}

\begin{figure}
    \centering
    \includegraphics[width=0.9\textwidth]{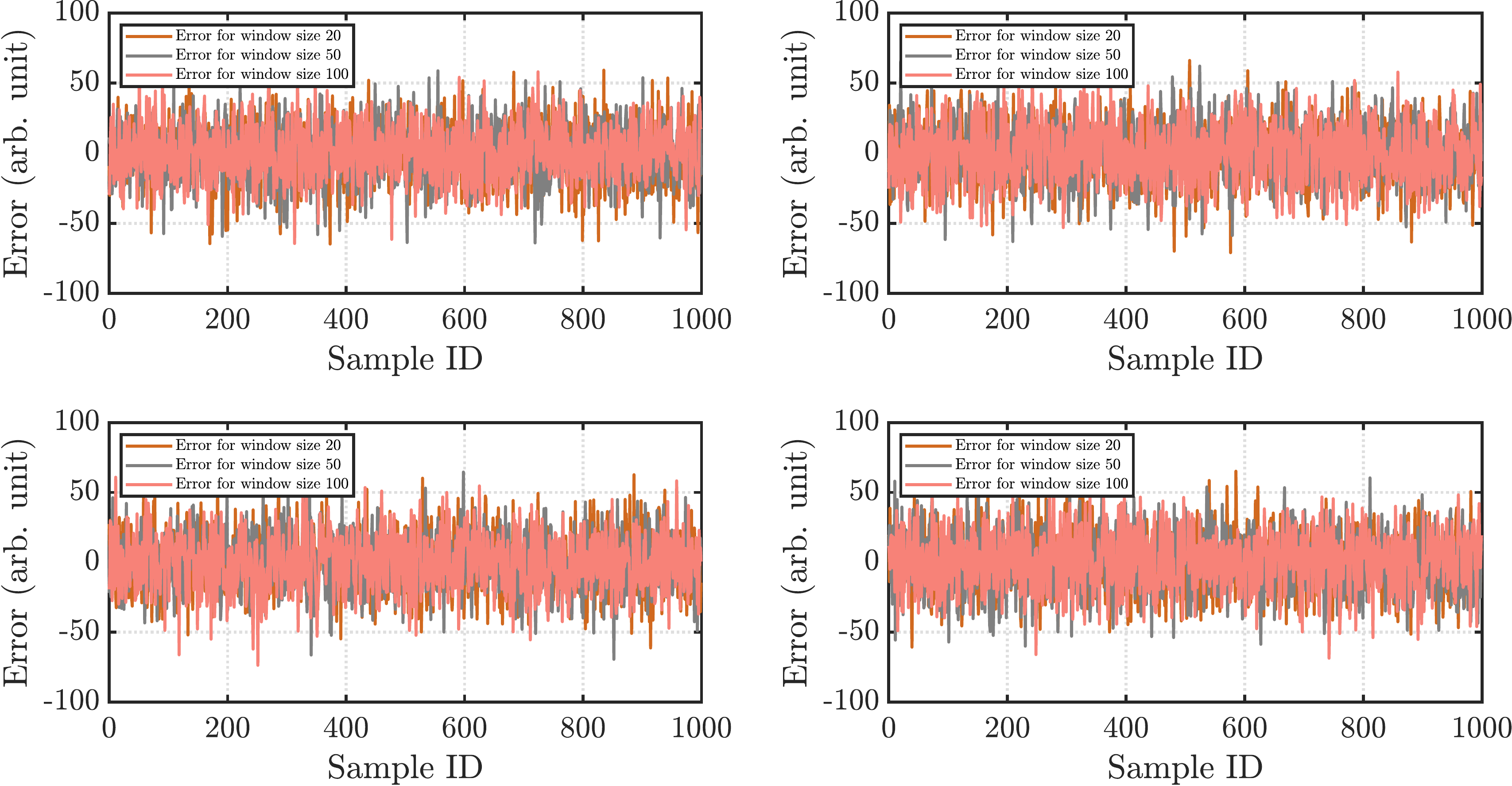}
    \caption{Variation of the error of the least squares prediction with respect to the observed data, with varying window sizes in the least squares optimization problems.}
    \label{fig:eeg_sysid_many}
\end{figure}

We present some preliminary results regarding the performance of the above approach. Specifically, we use $1000$ noisy measurements taken from $4$ channels of an intracranial electroencephalographic (iEEG) signal which records the brain activity of subjects undergoing epileptic seizures. The signals were recorded and digitized at a sampling rate of $512$ Hz at the Hospital of the University of Pennsylvania, Philadelphia, PA. Subdural grid and strip electrodes were placed at specialized locations (dictated by a multidisciplinary team of neurologists, neurosurgeons, and a radiologist), with the electrodes themselves consisting of linear and two-dimensional arrays spanning $2.3$ mm in diameter and having a inter-contact spacing of $10$ mm~\citep{khambhati2015dynamic,ashourvan2020model}. 

The least squares optimization problems described in Section~\ref{sec:system_identification} are solved using \texttt{CVX}~\citep{grant08,cvx} with the aid of a \mbox{window-based} approach using a finite subset of the entire range of measurements. This is done because the time series under consideration is nonlinear, and it is not possible to characterize the entire gamut of measurements using very few parameters. Figure~\ref{fig:eeg_sysid} shows the performance of our method on the above data. Additionally, we also show in Figure~\ref{fig:eeg_sysid_many} the variation of the error of the least squares predictions with respect to the observed data, with varying window sizes in the least squares optimization problems. We see that the identified system parameters are able to predict the system states fairly closely, thus demonstrating that our approach can be used to learn the system parameters of a discrete-time fractional-order system.

\subsection{Minimum-energy state estimation}

\begin{figure}[ht]
    \centering
    \includegraphics[width=0.9\textwidth]{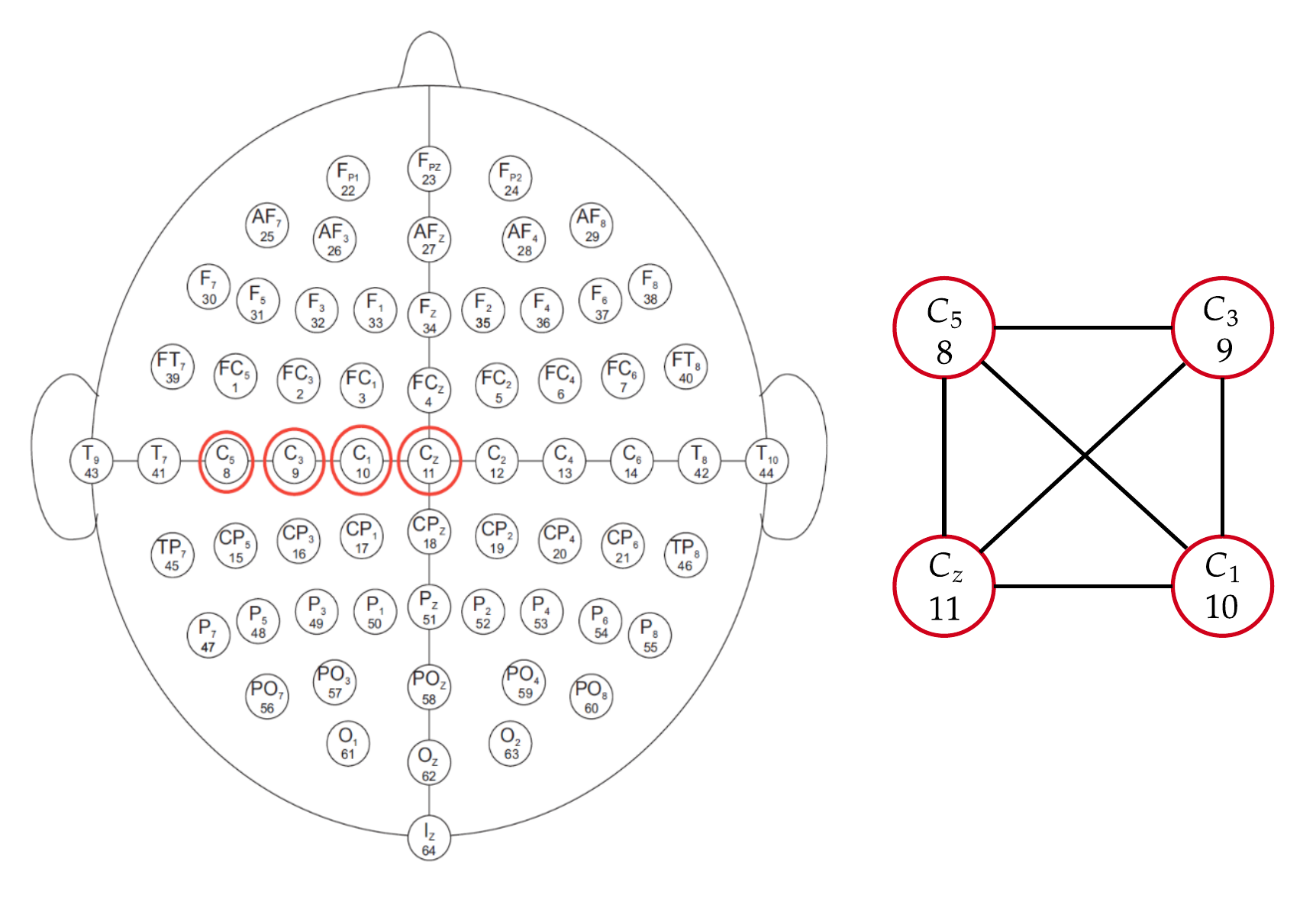}
    \caption{The distribution of the sensors for the measurement of EEG data is shown on the left. The channel labels are shown along with their corresponding numbers and the selected channels over the motor cortex are shown in red. The corresponding network formed by the EEG sensors is shown on the right.}
    \label{fig:eeg_networks}
\end{figure}

In this section, we consider the performance of the minimum-energy estimation paradigm on real-world neurophysiological networks considering EEG data. Specifically, we use $150$ noisy measurements taken from $4$ channels of a $64$-channel EEG signal which records the brain activity of subjects, as shown in Figure~\ref{fig:eeg_networks}. The subjects were asked to perform a variety of motor and imagery tasks. Furthermore, the specific choice of the $4$ channels was dictated due to them being positioned over the motor cortex of the brain, and, therefore, enabling us to predict motor actions such as the movement of the hands and feet. The data was collected using the BCI$2000$ system with a sampling rate of $160$ Hz~\citep{SchalkBCI,goldberger2000physiobank}. The spatial and temporal parameter components of the fractional-order system assumed to model the original EEG data were identified using the methods described in~\citep{gupta2018dealing}. The matrices $B_i = \begin{bmatrix}1&1&1&1 \end{bmatrix}^{\mathsf{T}}$ for all~$i$.

\begin{figure}[ht]
    \centering
    \includegraphics[width=0.9\textwidth]{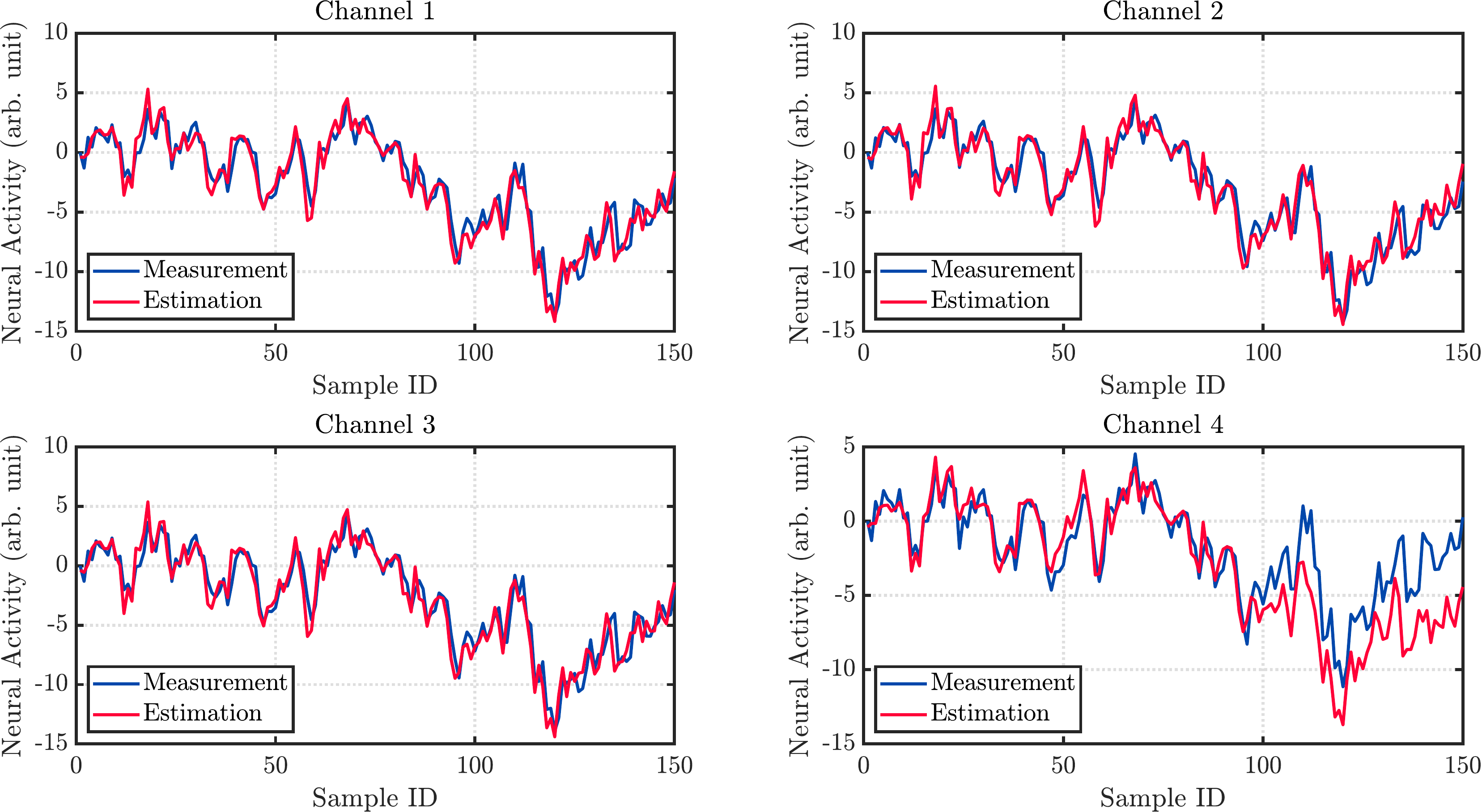}
    \caption{Comparison between the measured output of the $\mathfrak{v}$-augmented system (with $\mathfrak{v}=2$) versus the estimated output of a minimum-energy estimator implemented on the same, in the presence of process and measurement noises for $4$ channels of a $64$-channel EEG signal.}
    \label{fig:output_2}
\end{figure}

\begin{figure}[ht]
    \centering
    \includegraphics[width=0.9\textwidth]{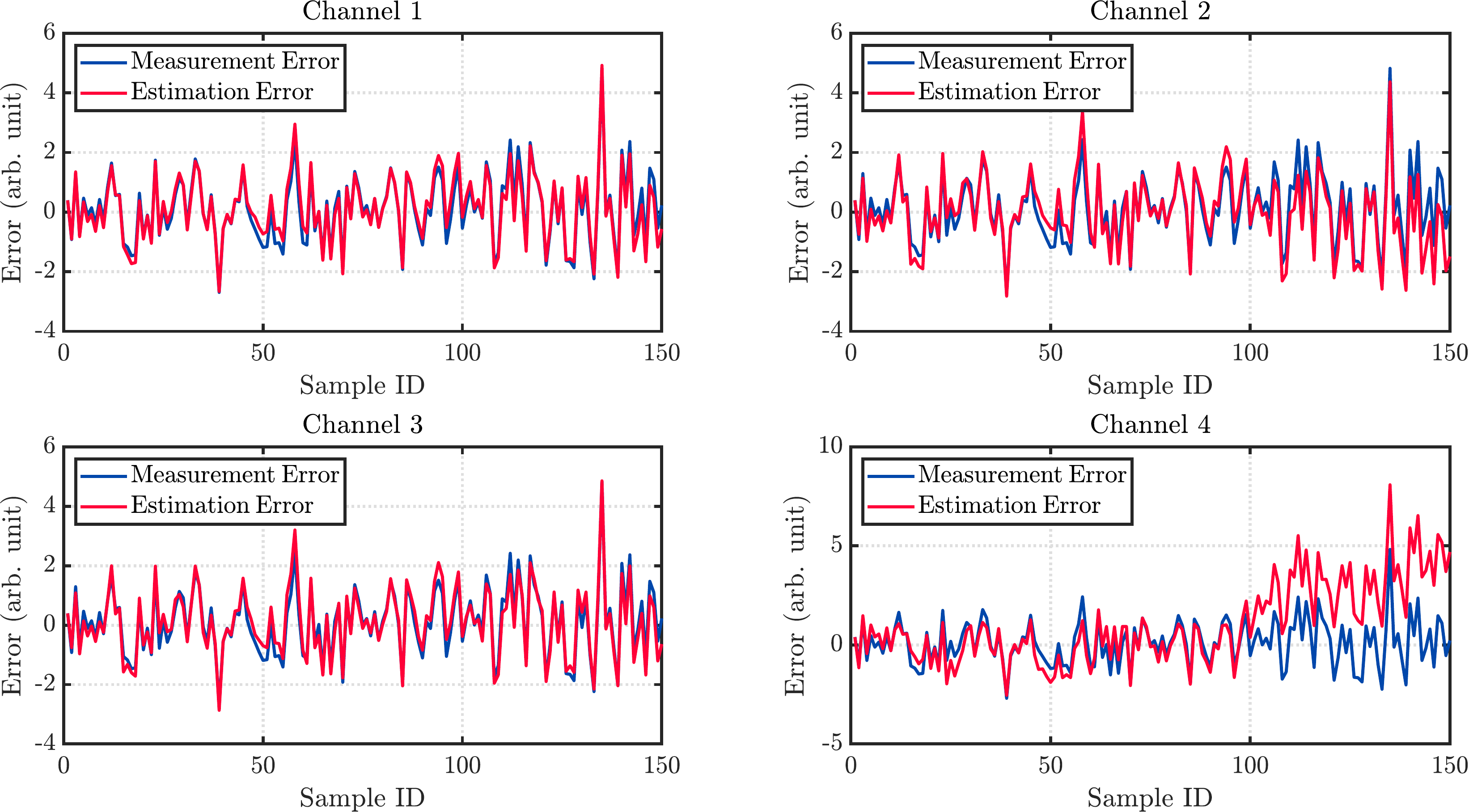}
    \caption{Comparison between the measurement error of the $\mathfrak{v}$-augmented system (with $\mathfrak{v}=2$) versus the estimation error of a minimum-energy estimator implemented on the same, in the presence of process and measurement noises for $4$ channels of a $64$-channel EEG signal.}
    \label{fig:error_2}
\end{figure}

\begin{figure}[ht]
    \centering
    \includegraphics[width=0.9\textwidth]{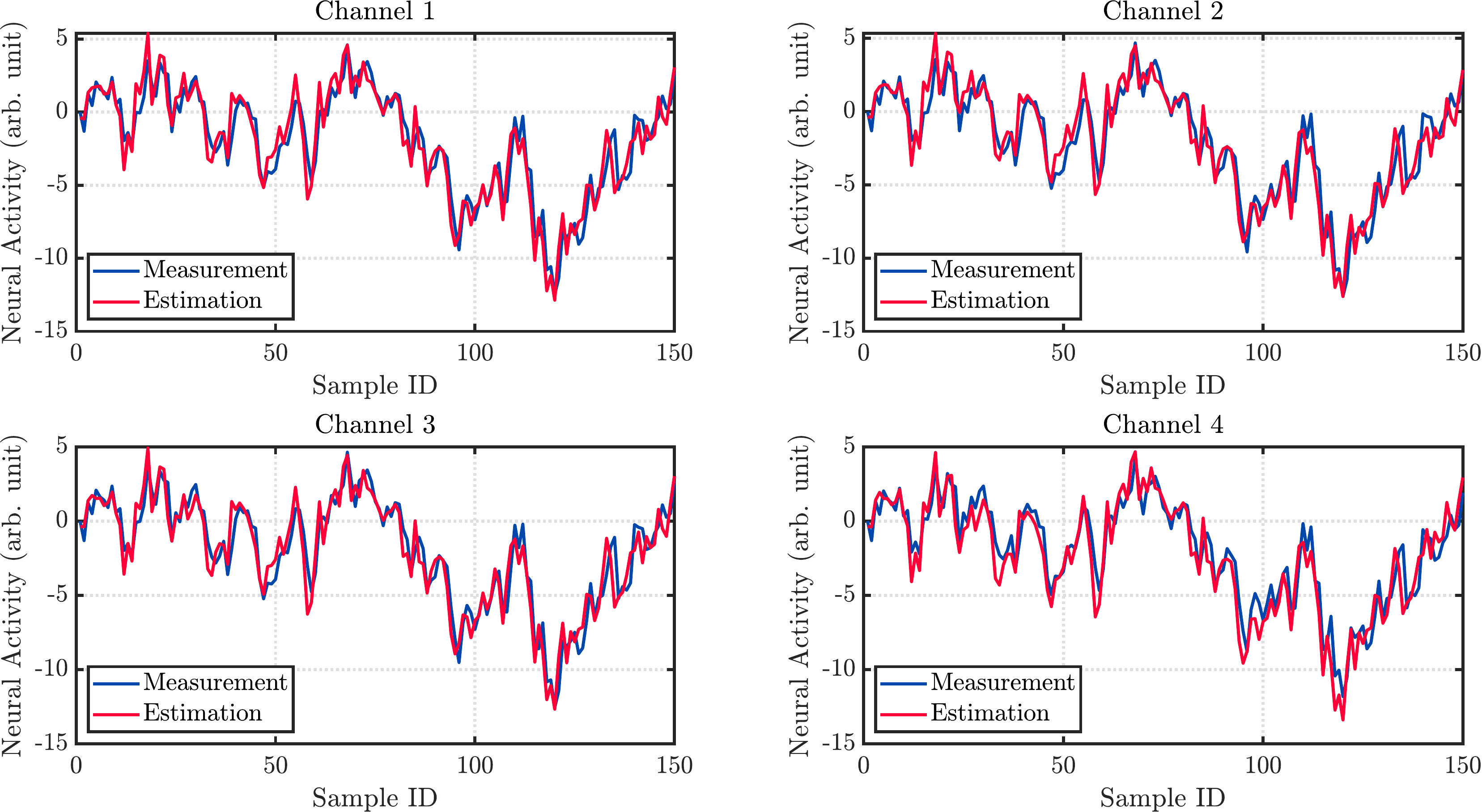}
    \caption{Comparison between the measured output of the $\mathfrak{v}$-augmented system (with $\mathfrak{v}=10$) versus the estimated output of a minimum-energy estimator implemented on the same, in the presence of process and measurement noises for $4$ channels of a $64$-channel EEG signal.}
    \label{fig:output_10}
\end{figure}

\begin{figure}[ht]
    \centering
    \includegraphics[width=0.9\textwidth]{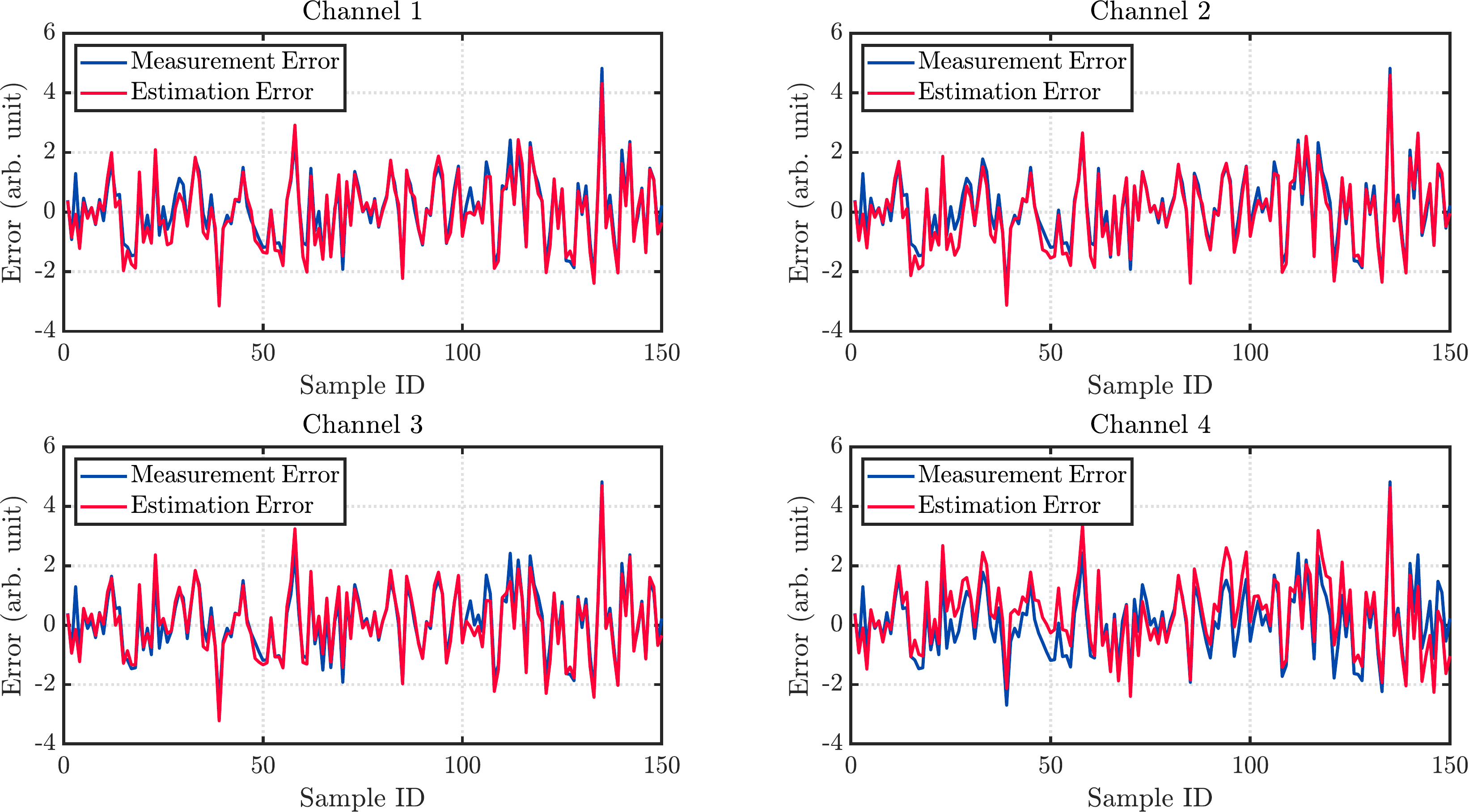}
    \caption{Comparison between the measurement error of the $\mathfrak{v}$-augmented system (with $\mathfrak{v}=10$) versus the estimation error of a minimum-energy estimator implemented on the same, in the presence of process and measurement noises for $4$ channels of a $64$-channel EEG signal.}
    \label{fig:error_10}
\end{figure}

\begin{figure}[ht]
    \centering
    \includegraphics[width=0.9\textwidth]{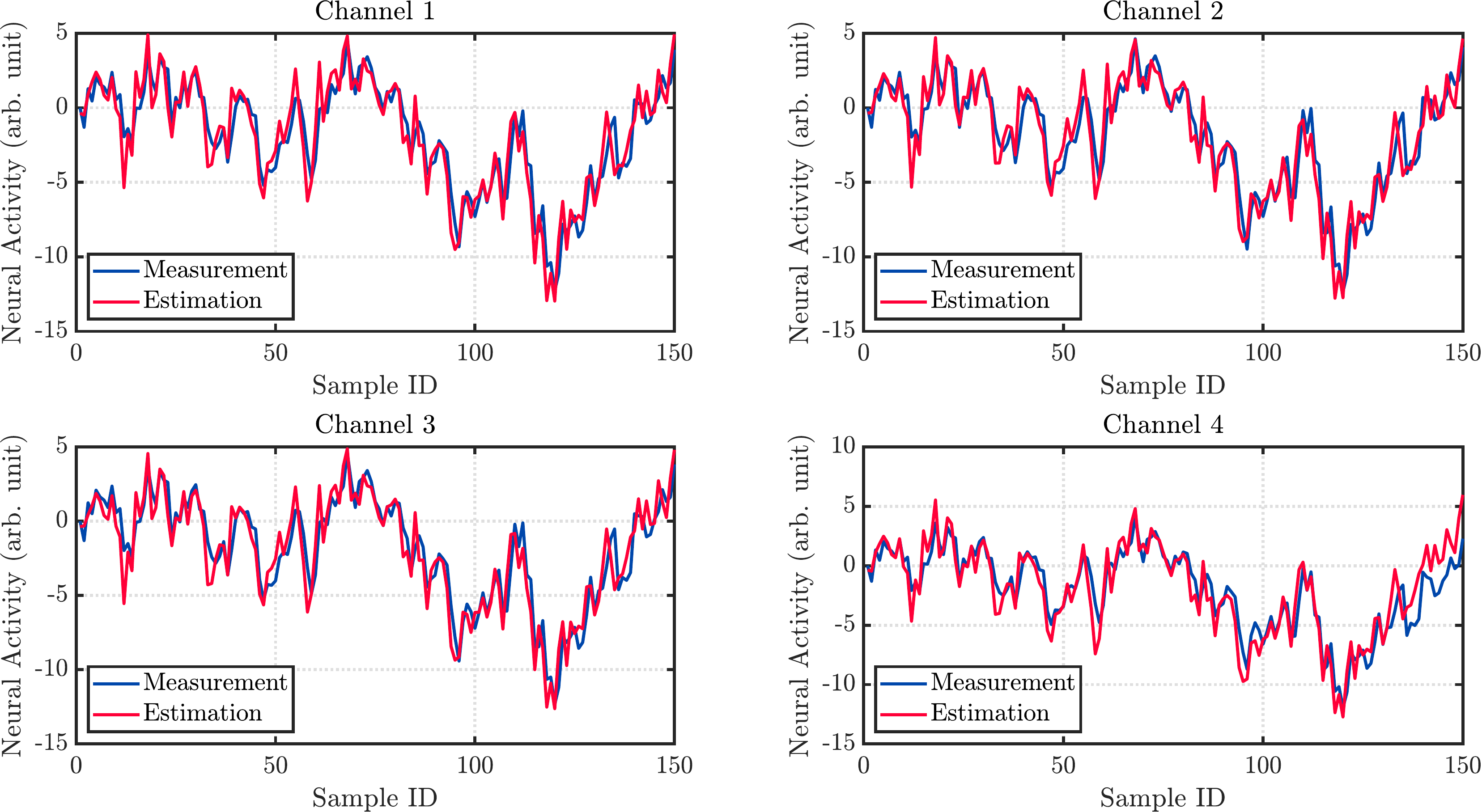}
    \caption{Comparison between the measured output of the $\mathfrak{v}$-augmented system (with $\mathfrak{v}=20$) versus the estimated output of a minimum-energy estimator implemented on the same, in the presence of process and measurement noises for $4$ channels of a $64$-channel EEG signal.}
    \label{fig:output_20}
\end{figure}

\begin{figure}[ht]
    \centering
    \includegraphics[width=0.9\textwidth]{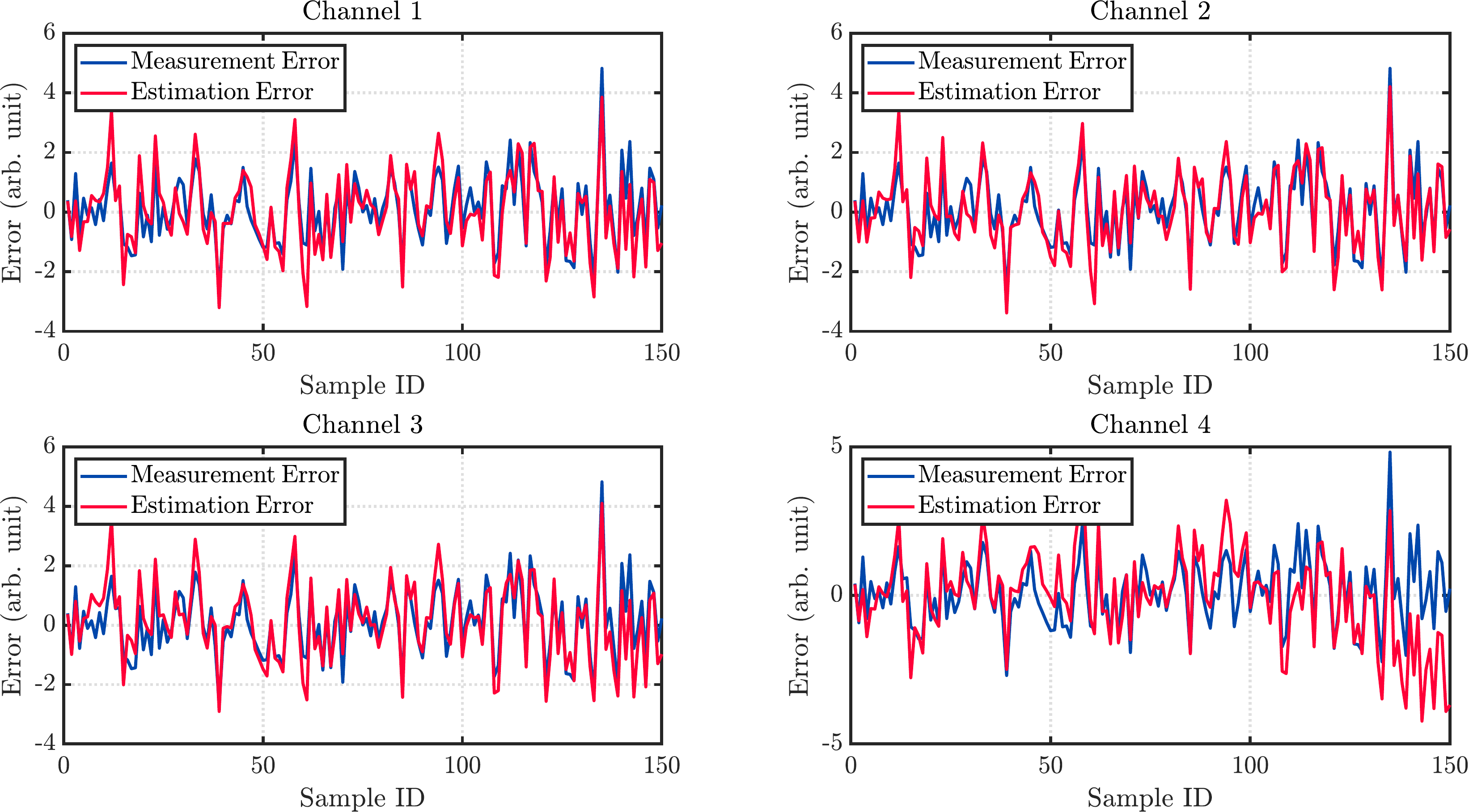}
    \caption{Comparison between the measurement error of the $\mathfrak{v}$-augmented system (with $\mathfrak{v}=20$) versus the estimation error of a minimum-energy estimator implemented on the same, in the presence of process and measurement noises for $4$ channels of a $64$-channel EEG signal.}
    \label{fig:error_20}
\end{figure}

The results of our approach, considering different values of $\mathfrak{v}$, are shown in Figures~\ref{fig:output_2} and~\ref{fig:error_2} (for $\mathfrak{v}=2$), Figures~\ref{fig:output_10} and~\ref{fig:error_10} (for $\mathfrak{v}=10$), and Figures~\ref{fig:output_20} and~\ref{fig:error_20} (for $\mathfrak{v}=20$), which show, respectively (for each value of $\mathfrak{v}$), the comparison between the measured output of the network with noise and the estimated response obtained from the \mbox{minimum-energy} estimator, and also the juxtaposition of the measurement error and the estimation error of the \mbox{minimum-energy} estimation process. We find that the minimum-energy estimator is successfully able to estimate the states in the presence of noise in both the dynamics and the measurement processes.

We also note from the Figures~\ref{fig:output_2} and~\ref{fig:error_2} that when $\mathfrak{v}=2$, we get comparatively larger estimation errors associated with the last $50$ or so samples of Channel $4$ and that this behavior can be mitigated by increasing the value of $\mathfrak{v}$, e.g., by choosing $\mathfrak{v}=10$ or $\mathfrak{v}=20$. This is in line with the discussion at the end of Section~\ref{sec:state_estimation}, and choosing a larger value of $\mathfrak{v}$ can always, in practice, provide us with better estimation performances, as seen from this example. 

\subsection{Neurostimulation using fractional-order model predictive control for epileptic seizure mitigation}

In what follows, we propose to illustrate the use of the fractional-order system model predictive control (FOS-MPC) framework for neurostimulation in the context of mitigating epileptic seizures. We demonstrate the workings of the proposed approach on four different experimental scenarios relying primarily on intracranial electroencephalographic (iEEG) data: 
\emph{(i)}~an iEEG signal demonstrating an epileptic seizure simulated by the neural mass model proposed by Jansen and Rit~\mbox{\citep{jansen1993neurophysiologically,jansen1995electroencephalogram}}; 
\emph{(ii)}~an iEEG signal simulated by a neural field model proposed by Martinet et al. in~\citep{martinet2017human} that replicates the spatiotemporal dynamics of a seizure; 
\emph{(iii)}~an iEEG signal simulated by the phenomenological `Epileptor' model proposed in~\citep{jirsa2014nature}; and 
\emph{(iv)}~real-time iEEG signals for three human subjects undergoing epileptic seizures. 
For all of the above cases, we start by considering an epileptic seizure, captured by a linear fractional-order system (FOS) model, whose parameters are obtained through a system identification method using brainwave data obtained from iEEGs. 

\subsubsection{Epileptic seizure simulated by the Jansen-Rit neural mass model}

Although initially proposed to account for human EEG rhythms and visual evoked potentials, the Jansen-Rit neural mass model~\citep{jansen1995electroencephalogram} has also been used to shed light on human epileptiform brain dynamics~\mbox{\citep{wendling2000relevance,wendling2016computational}.} The Jansen-Rit neural mass model is composed of three interacting subpopulations that include: 
the main subpopulation, the excitatory feedback subpopulation, and the inhibitory feedback subpopulation. The structure of the model is such that the main subpopulation comprises cells that receive neuronal signals in feedback from the excitatory and inhibitory subpopulations. 


The use of neural mass models akin to the Jansen-Rit model in feedback control frameworks is well documented. All the works in~\citep{wang2016suppressing,xia2019closed,wei2019seizure,wei2019seizure_2,wei2019control,soltan2018fractional} use neural mass models, in the control theory sense, for the suppression of epileptic seizures. In what follows, we will demonstrate the effectiveness of our proposed control strategy on a seizure simulated by the classical Jansen-Rit neural mass model with standard parameter values.

First, we need to determine the  parameters $A$ and $\alpha$ that model both spatial coupling and fractional coefficients, respectively, that craft the evolution of the state $x[k] \in\mathbb{R}^n$ in the fractional-order system (FOS) model.
\begin{equation}
    \Delta^\alpha x[k+1] = Ax[k] + B u[k] + B^w w[k],
    \label{eq:sysIDed_jne}
\end{equation}
with $w[k]$ denoting additive white Gaussian noise (AWGN). Since the system is single-input-single-output (SISO), we have both $A$ and $\alpha$ to be scalars. To identify the parameters $A$ and $\alpha$, we used the method proposed in~\citep{gupta2018dealing}. The parameters obtained are $A = -0.0054$ and $\alpha = 1.4881$. Furthermore, we assume that $B = 1$ and $B^w = 0.1$.

For the cost function in \eqref{eq:mpc123_jne}, we utilized $Q_k=I_n$, $R_k=I_{n_u}$, and $c_k=0_{n_u\times 1}$ (with $n = n_u = 1$), to emphasize minimizing the overall energy in the measured iEEG signal, while penalizing slightly for overly aggressive stimulation. Furthermore, we included a safety linear constraint of $-5 \leq u[k] \leq 5$. Our predictive model was based on a $(p=15)-$step (15 ms) predictive model approximation of the FOS plant, with a $(P=20)-$step (20 ms) prediction horizon and $(M=10)-$step (10 ms) control horizon.
\begin{figure}
    \centering
    \includegraphics[width=0.9\textwidth]{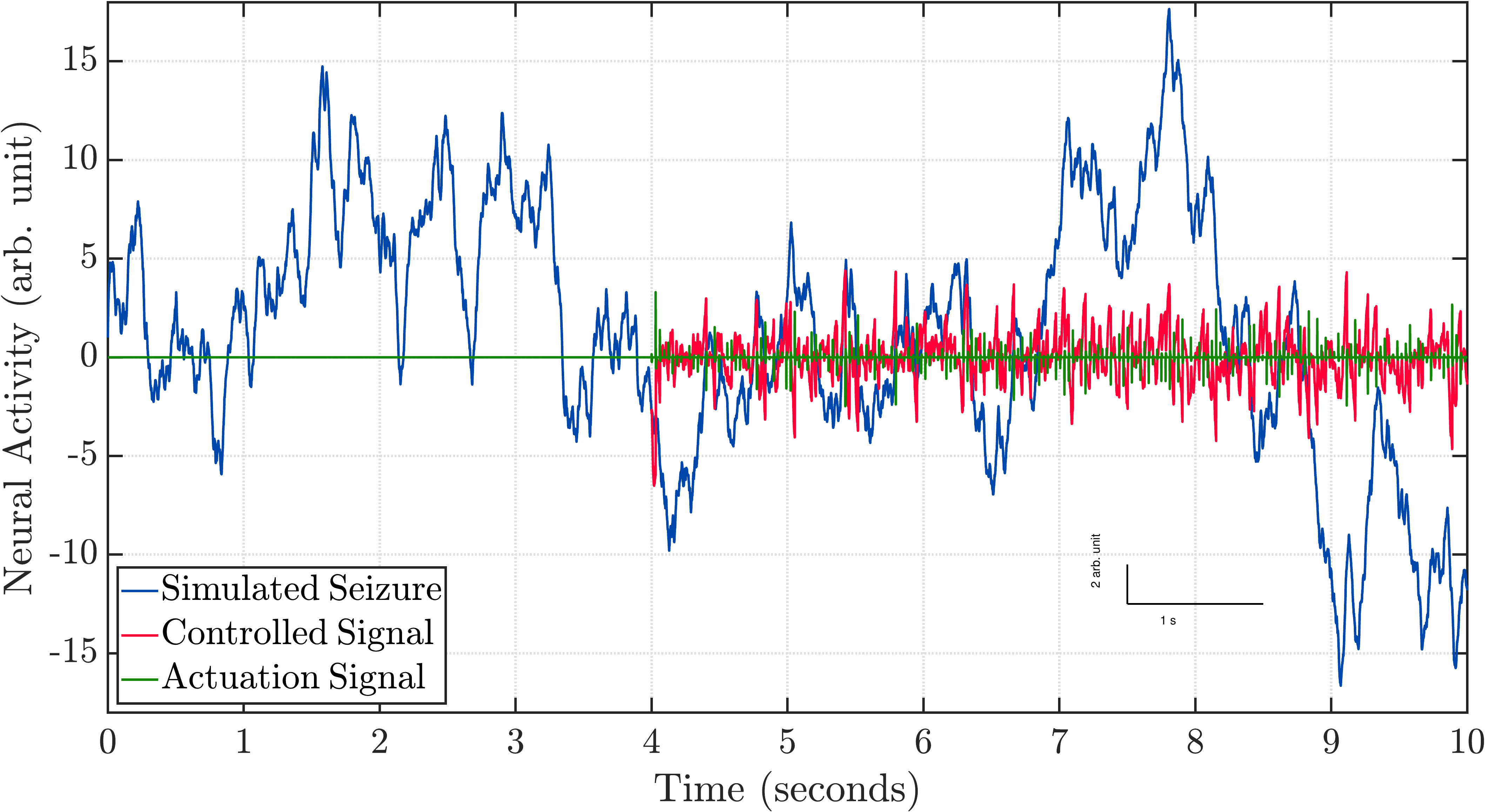}
    \caption{Results of the proposed FOS-MPC closed-loop neurostimulation strategy on an iEEG seizure simulated by the Jansen-Rit neural mass model. The simulated iEEG signal with the seizure is depicted in blue, the controlled signal is depicted in red, and the stimulation pulses are shown in green.}
    \label{fig:nmm_mpc}
\end{figure}
The results are presented in Figure~\ref{fig:nmm_mpc}, which provide evidence that the proposed stimulation strategy allows us to achieve amplitude suppression using a (time-varying) impulse-like stimulation scheme. Note that the actuation signal $u_k$ kicks in at about the 4-second mark in the figure.

\subsubsection{Epileptic seizure simulated by the mean-field model proposed by Martinet et al.~\citep{martinet2017human}.}

Next, we turn our attention towards a computational model that uses traveling wave dynamics to capture inter-scale coupling phenomena between large-scale neural populations in the cortex and small-scale groups in cortical columns~\citep{martinet2017human}. Modeling the complex spatiotemporal dynamics of epileptic seizures is a challenging task, mainly because of the interaction of myriad scales in both time and space.

The neural field model proposed by Martinet et al. in~\citep{martinet2017human} is a modified version of the mean-field model proposed in~\citep{steyn2013interacting} that seeks to explain the phenomena, origin, and spatiotemporal dynamical properties of seizure propagation and spike-and-wave discharges (SWDs). Additionally, their work advances the hypothesis that increased diffusion of extracellular potassium concentrations in space influences the interlaced coupling of human seizures. In what follows, we will use the simulated seizure data obtained from the aforementioned model and then consider our closed-loop MPC neurostimulation scheme on the same model.

\begin{figure}
    \centering
    \includegraphics[width=0.9\textwidth]{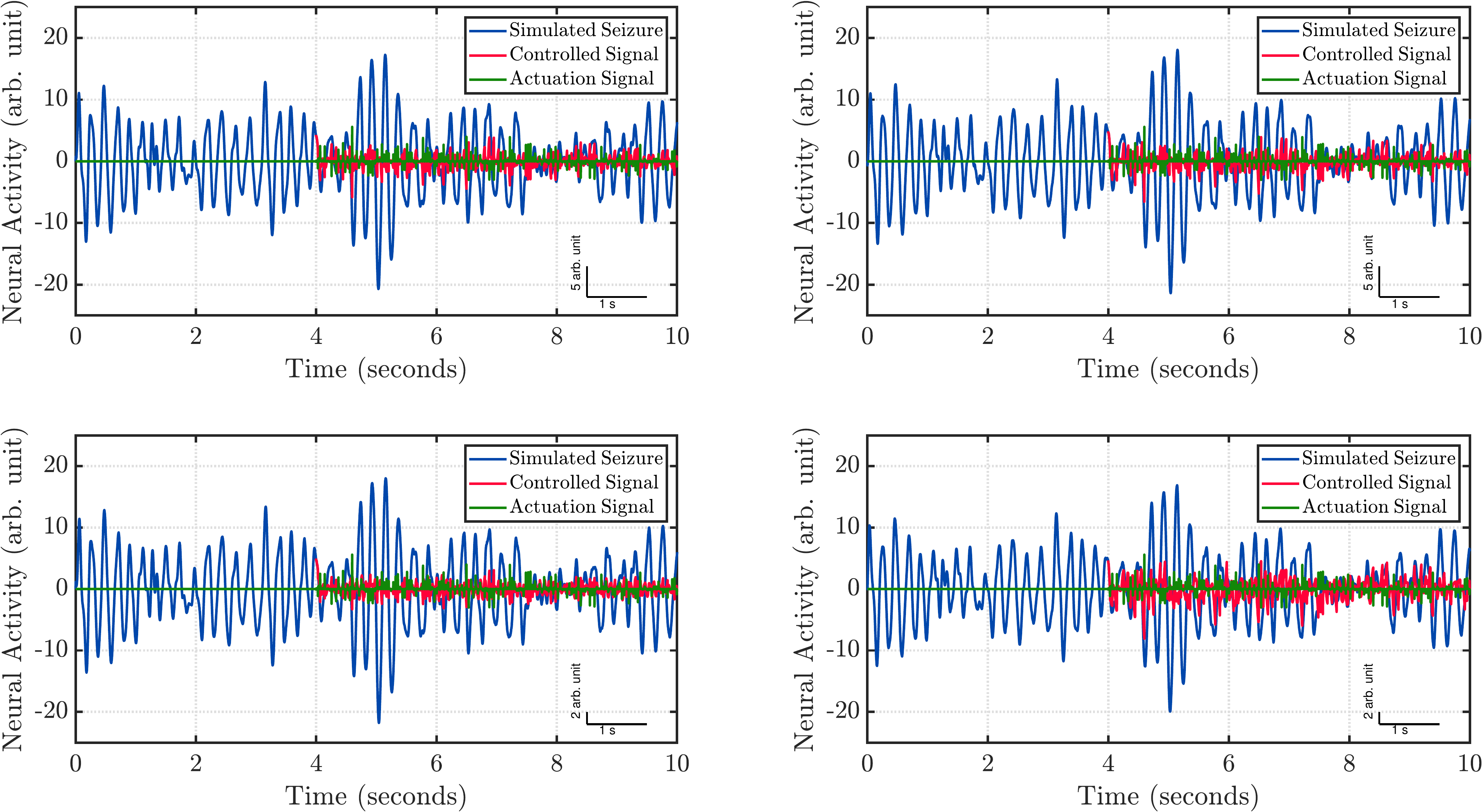}
    \caption{Results of the proposed FOS-MPC closed-loop neurostimulation strategy on an iEEG seizure simulated by the traveling wave dynamics model proposed in Martinet et al.~\citep{martinet2017human}. The simulated iEEG signal with the seizure is depicted in blue, the controlled signal is depicted in red, and the stimulation pulses are shown in green.}
    \label{fig:nmm_kramer_mpc}
\end{figure}

To determine the system parameters $A$ and $\alpha$ in~\eqref{eq:sysIDed_jne}, we utilize roughly 2 seconds of pre-ictal activity captured by the model. Note that here, we will only consider $n=4$ channels for our proposed approach to mimic the capabilities available in the NeuroPace\textsuperscript{\textregistered} RNS\textsuperscript{\textregistered} device. Applying the methods in~\citep{gupta2018dealing} yields the following FOS parameters:
\begin{equation}
    A = \begin{bmatrix}
    0.2969 &  -0.0203 &  -0.2922 & 0.0587\\
    0.2574  & -0.1726 &  -0.1905 &    0.1535\\
    0.5348  & -0.1066 &  -0.3471  & -0.0169\\
    0.4007  & -0.6752 &   0.0044  &  0.3186     
    \end{bmatrix},
\end{equation}
and
\begin{equation}
    \alpha = \begin{bmatrix} 0.8114 & 0.8334 & 0.8034 & 0.8413 \end{bmatrix}^{\mathsf{T}}.
\end{equation}
Additionally, we consider a single control signal $u_k$ that affects all the channels equally, i.e., $B = \begin{bmatrix} 1&1&1&1 \end{bmatrix}^{\mathsf{T}}$ and the matrix of weights $B^w = 0.05 I_4$, with $I_4$ being the $4 \times 4$ identity matrix.

Using the FOS-MPC neurostimulation strategy with $Q_k=I_n$, $R_k=I_{n_u}$, and $c_k=0_{n_u\times 1}$ (with $n = 4$ and $n_u = 1$), and safety linear constraints of $-100 \leq u[k] \leq 100$, we find from Figure~\ref{fig:nmm_kramer_mpc} that our proposed approach successfully suppresses seizure-like activity using a (time-varying) impulse-like stimulation scheme. In this case, we use a $(p=10)-$step (20 ms) predictive model approximation of the FOS plant, with a $(P=10)-$step (20 ms) prediction horizon, and $(M=8)-$step (16 ms) control horizon. Here too, the actuation signal $u_k$ kicks in at about the 4-second mark.

\subsubsection{Epileptic seizure simulated by the Epileptor, a phenomenological model of seizures by Jirsa et al.~\citep{jirsa2014nature}}

Next, we investigate the performance of our proposed approach on the Epileptor model~\citep{jirsa2014nature}, which is a phenomenological model able to accurately reproduce the dynamics of a wide variety of human epileptic seizures recorded with iEEG electrodes.

The Epileptor is a mathematical model proposed by Jirsa et al. in~\citep{jirsa2014nature} and is based on analyzing experimental readings of iEEG seizure discharges in various human and animal subjects. At its core, the model consists of six coupled ordinary differential equations in three time scales. These equations are successfully able to model bistable dynamics between alternating fast discharges and inter-ictal activity, spike-and-wave events (SWEs), and the evolution of the neural populations through the phenomena of seizure onset and offset. In what follows, we will use the simulated seizure data obtained from the Epileptor model and implement our closed-loop MPC neuromodulation scheme on it.

To determine the parameters $A$ and $\alpha$ that model both spatial coupling and fractional coefficients, respectively, that craft the evolution of the state dynamics in~\eqref{eq:sysIDed_jne}, we use the method proposed in~\citep{gupta2018dealing}. Here, like the Jansen-Rit model, the system is SISO, and hence $A$ and $\alpha$ are scalars. The parameters obtained are $A = -0.0051$ and $\alpha = 1.0614$. Furthermore, we assume that $B = 1$ and $B^w = 0.25$.

We implement the FOS-MPC neurostimulation strategy with $Q_k=I_n$, $R_k=I_{n_u}$, and $c_k=0_{n_u\times 1}$ (with $n = n_u = 1$) and safety linear constraints of $-50 \leq u[k] \leq 50$. In this case, our predictive model was based on a $(p=20)-$step predictive model approximation of the FOS plant, with a $(P=20)-$step prediction horizon and $(M=10)-$step control horizon.
\begin{figure}
    \centering
    \includegraphics[width=0.9\textwidth]{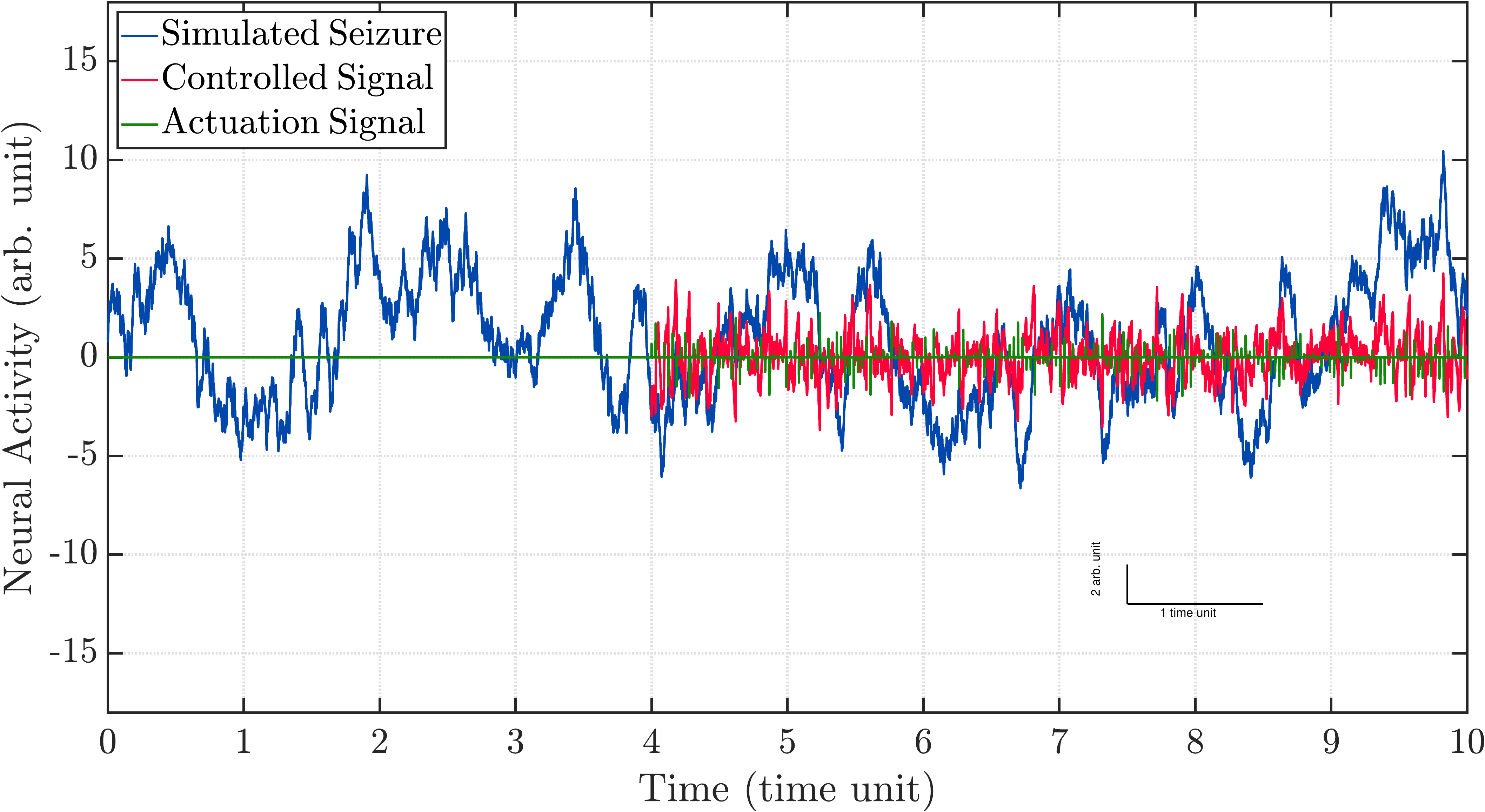}
    \caption{Results of the proposed FOS-MPC closed-loop neurostimulation strategy on an iEEG seizure simulated by the Epileptor model. The simulated iEEG signal with the seizure is depicted in blue, the controlled signal is depicted in red, and the stimulation pulses are shown in green.}
    \label{fig:epileptor_simu}
\end{figure}
The results are presented in Figure~\ref{fig:epileptor_simu}, which provide evidence that the proposed stimulation strategy allows us to achieve amplitude suppression for a seizure simulated by the Epileptor model with standard parameter values.

\section{Conclusions and directions for future research}\label{sec:future_research}

Cyber-neural systems are becoming pervasive in today's society, yet they still lack the capability of performing real-time closed-loop control of neural activity. Control systems engineers will play a vital role in bringing this technology to reality as they develop the tools required by interdisciplinary teams involved in envisioning the next generation of neurotechnology.

That said, we provided a glance at some of the latest trends and techniques in \mbox{fractional-order} based system modeling, analysis, and closed-loop control towards the development of future neurotechnologies. In particular, we present results on system identification,  state estimation, and closed-loop control in discrete-time fractional-order dynamical systems.

Notwithstanding, there is a plethora of interesting research directions which can be taken from here, aligned with each of the foundations pinpointed next.  

\subsection{System identification}

System identification of \mbox{fractional-order} systems is an extremely \mbox{under-explored} field in general, with a lack of a systematic and unified theory, with some preliminary approaches utilizing wavelets~\citep{flandrin1992wavelet}, \mbox{frequency-domain} techniques~\citep{adams2006fractional,dzielinski2011identification}, or a sequential combination of wavelets and expectation-maximization (EM)~\citep{gupta2018dealing}.

 Although approaches using the EM algorithm for linear~\citep{gibson2005robust} as well as nonlinear system identification~\citep{schon2011system} have existed in the literature for a while now, one immediately notices that there is a \mbox{long-standing} problem in characterizing theoretical robustness guarantees for the same.
We can draw inspiration from preliminary analyses of \mbox{finite-sample} robustness guarantees for EM in~\citep{wu2016convergence,balakrishnan2017statistical,yan2017convergence} to characterize the sample complexity in identifying linear \mbox{time-invariant} (LTI) to do the same for the \mbox{fractional-order} systems.

In practice, it would be important to investigate approaches to identify the spatial and temporal parameters of fractional-order systems based on bootstrapping~\citep{tjarnstrom1999use}, with alternating and progressively better identifications of the aforementioned parameters. Future work should focus on chalking out a general theory of identifying certain classes of \mbox{fractional-order} systems.

Alternatively, one could also look into strategies behind using recurrent neural networks (RNNs) to identify the fractional-order systems' parameters. One of the most celebrated results in neural network theory, the \emph{universal approximation theorem}~\citep{hornik1989multilayer,funahashi1989approximate,cybenko1989approximation}, states that continuous functions can be arbitrarily well-approximated by single-hidden-layer feedforward neural networks. While recent work~\citep{hutter2021metric} seems to suggest the presence of allied results when RNNs are used to identify stable LTI systems optimally in the sense of metric entropy~\citep{zames1993note}, it remains to be seen whether universal approximation theorem-like results can be derived when RNNs are used to identify fractional-order systems.

It would also be interesting to investigate fundamental \mbox{information-theoretic} connections between the number of samples needed to perform online system identification for fractional-order systems that would, in turn, allow for a more robust control design using the same. Furthermore, one could also potentially look into the number of samples needed to uniquely identify the parameters of fractional-order systems and whether different identified realizations potentially correspond to different fractional-order systems. 

\subsection{State estimation}

Although fractional-order systems have found vast success in modeling the spatiotemporal properties of EEG, some of the properties accounted for by these models actually originate from unknown sources external to the system under consideration. 
Future work should focus on modeling these external sources by unknown input stimuli and then focus on state estimation of the resultant model with unknown inputs. 

Real-time neural activity can be monitored to \mbox{self-regulate} brain function. This is known in the literature as neurofeedback~\citep{marzbani2016neurofeedback}, and it would be interesting to study how the introduction of feedback to such a system changes our perspectives on this problem.

Besides, whereas the construction of resilient state estimators grow over the last decade, little effort was put in developing resilient versions of state estimators for fractional-order systems. In particular, and given that suitable assumptions for the disturbance and noise do not rely on gaussianity assumptions, it would be imperative to build a resilient and attack-resistant version of the minimum-energy estimator. Specifically, to consider adversarial attacks or artifacts associated with the measurement process, since the former approach is consistent with the fact that (adversarial) attacks on sensors often do not follow any particular dynamic or stochastic characterization. 

Last but not least, it would be crucial to understand how to design filter-like approaches that amalgamate the problems of simultaneous system identification and estimation suitable for the deployment in real-time CNS.

\subsection{Closed-loop control}

Very rarely in practical settings do we have deterministic \mbox{fractional-order} models. As we saw, neural signals are particularly prone to artifacts from outside the brain. Furthermore, stabilizing the underlying models in the presence of disturbances becomes relevant in the treatment of disorders like epilepsy, Parkinson's, or Alzheimer's disease. 

In recent years, there have been increasing research efforts into finding possible therapies for the aforementioned using neurofeedback~\citep{marzbani2016neurofeedback}. Future work, therefore, should focus on developing controllers and observers for fractional-order systems with the associated process and measurement noise and investigating the possible existence of separation \mbox{principle-like} results akin to those already existing in the field of linear stochastic control theory.

Another direction of work entails deriving robustness guarantees for controlling \mbox{discrete-time} fractional-order systems using an \mbox{inner-outer} loop control strategy. Specifically, in this context, we seek to discover the advantages and disadvantages of truncating a \mbox{discrete-time} fractional-order systems according to a given truncation horizon, thus approximating the fractional-order systems as an augmented LTI system and performing model predictive control with the same. 

Additionally, one could also rely on some tools from robust control, namely integral quadratic constraints (IQCs)~\citep{megretski1997system}. 
IQCs are, essentially, inequalities used to describe possible input-output signals resulting from a system component that is challenging to model because it is either nonlinear, time-varying, noisy, or switch stochastically or adversarially with time. A particular issue of interest will be to explore the \mbox{trade-offs} in performance when the fractional-order systems (which represents the inner loop) are written as an augmented LTI system due to a fixed truncation horizon versus when it is modeled as a non-Markovian nonlinearity with IQCs characterizing the same.

Additionally, although finite-time LTI truncations of fractional-order systems with constant truncation horizons are considered in this paper, fractional-order systems inherently possess infinite long-term memory. The question that is an immediate consequence of the latter fact is whether the theory of linear control systems in infinite dimensions~\citep{curtain2012introduction} can be used to provide key insights into control-theoretic properties such as controllability, observability, and stabilizability for such systems. While there have been some preliminary works in this direction~\citep{baleanu2019approximate,zitane2020stability,wei2019analysis,sabatier2021fractional}.
Consequently, future work must consist of using mathematical techniques used in the analysis of partial and delay differential equations, in particular, operator equations and $C_0$-semigroup theory~\citep{bamieh2002distributed} for fractional-order systems.

From an engineering or applied control point-of-view, it is important to pinpoint some limitations and drawbacks of current \mbox{event-triggered} open-loop stimulation strategies (i.e., they can be inefficient or even cause seizure-like activity). Consequently, it serves as a call for action from neurophysiologists and engineers that work with neurostimulation  (as well as deep brain stimulation) devices, towards validation in \emph{in vitro} and \emph{in vivo} scenarios.  %
That said,  the advances in computational processing power made in the last 10-20 years have made the prospects of turning into reality technology that was theoretically devised and previously impossible to implement in real-life. MPC and fractional-order systems-based technologies both fall under this category and have thus been significantly overlooked in the industry. However, both are growing in popularity amongst several research communities, and some predict a considerably more widespread impact than originally thought.  

Notwithstanding, the validation is insufficient to establish a framework since several foundational problems need to be addressed. Specifically, the robustness of the stimulation strategies concerning the parameters of the models (e.g., the dynamics and the stimuli deployed, as well as the approximations considered to attain real-time stimulation) in devices with low storage, and limited battery and computation capabilities. Towards this goal, only transdisciplinary work between scientists and engineers will lead to success that ultimately will be reflected in the quality of life improvement of patients with neurological disorders (e.g., epilepsy).



\bibliographystyle{elsarticle-harv.bst} 
\bibliography{all_phd_refs}





\end{document}